\documentclass{siamart1116}

\newcommand{\rulesep}{\unskip\ \vrule\ }
\renewcommand{\emph}[1]{\textsl{#1}}

\usepackage{lipsum}
\usepackage{amsfonts}
\usepackage{graphicx}
\usepackage{subcaption}
\usepackage{epstopdf}
\usepackage{csquotes}
\usepackage{xfrac}
\usepackage{listings}
\usepackage[draft]{minted}

\makeatletter
\def\PYGdefault@reset{\let\PYGdefault@it=\relax \let\PYGdefault@bf=\relax%
    \let\PYGdefault@ul=\relax \let\PYGdefault@tc=\relax%
    \let\PYGdefault@bc=\relax \let\PYGdefault@ff=\relax}
\def\PYGdefault@tok#1{\csname PYGdefault@tok@#1\endcsname}
\def\PYGdefault@toks#1+{\ifx\relax#1\empty\else%
    \PYGdefault@tok{#1}\expandafter\PYGdefault@toks\fi}
\def\PYGdefault@do#1{\PYGdefault@bc{\PYGdefault@tc{\PYGdefault@ul{%
    \PYGdefault@it{\PYGdefault@bf{\PYGdefault@ff{#1}}}}}}}
\def\PYGdefault#1#2{\PYGdefault@reset\PYGdefault@toks#1+\relax+\PYGdefault@do{#2}}

\expandafter\def\csname PYGdefault@tok@w\endcsname{\def\PYGdefault@tc##1{\textcolor[rgb]{0.73,0.73,0.73}{##1}}}
\expandafter\def\csname PYGdefault@tok@c\endcsname{\let\PYGdefault@it=\textit\def\PYGdefault@tc##1{\textcolor[rgb]{0.25,0.50,0.50}{##1}}}
\expandafter\def\csname PYGdefault@tok@cp\endcsname{\def\PYGdefault@tc##1{\textcolor[rgb]{0.74,0.48,0.00}{##1}}}
\expandafter\def\csname PYGdefault@tok@k\endcsname{\let\PYGdefault@bf=\textbf\def\PYGdefault@tc##1{\textcolor[rgb]{0.00,0.50,0.00}{##1}}}
\expandafter\def\csname PYGdefault@tok@kp\endcsname{\def\PYGdefault@tc##1{\textcolor[rgb]{0.00,0.50,0.00}{##1}}}
\expandafter\def\csname PYGdefault@tok@kt\endcsname{\def\PYGdefault@tc##1{\textcolor[rgb]{0.69,0.00,0.25}{##1}}}
\expandafter\def\csname PYGdefault@tok@o\endcsname{\def\PYGdefault@tc##1{\textcolor[rgb]{0.40,0.40,0.40}{##1}}}
\expandafter\def\csname PYGdefault@tok@ow\endcsname{\let\PYGdefault@bf=\textbf\def\PYGdefault@tc##1{\textcolor[rgb]{0.67,0.13,1.00}{##1}}}
\expandafter\def\csname PYGdefault@tok@nb\endcsname{\def\PYGdefault@tc##1{\textcolor[rgb]{0.00,0.50,0.00}{##1}}}
\expandafter\def\csname PYGdefault@tok@nf\endcsname{\def\PYGdefault@tc##1{\textcolor[rgb]{0.00,0.00,1.00}{##1}}}
\expandafter\def\csname PYGdefault@tok@nc\endcsname{\let\PYGdefault@bf=\textbf\def\PYGdefault@tc##1{\textcolor[rgb]{0.00,0.00,1.00}{##1}}}
\expandafter\def\csname PYGdefault@tok@nn\endcsname{\let\PYGdefault@bf=\textbf\def\PYGdefault@tc##1{\textcolor[rgb]{0.00,0.00,1.00}{##1}}}
\expandafter\def\csname PYGdefault@tok@ne\endcsname{\let\PYGdefault@bf=\textbf\def\PYGdefault@tc##1{\textcolor[rgb]{0.82,0.25,0.23}{##1}}}
\expandafter\def\csname PYGdefault@tok@nv\endcsname{\def\PYGdefault@tc##1{\textcolor[rgb]{0.10,0.09,0.49}{##1}}}
\expandafter\def\csname PYGdefault@tok@no\endcsname{\def\PYGdefault@tc##1{\textcolor[rgb]{0.53,0.00,0.00}{##1}}}
\expandafter\def\csname PYGdefault@tok@nl\endcsname{\def\PYGdefault@tc##1{\textcolor[rgb]{0.63,0.63,0.00}{##1}}}
\expandafter\def\csname PYGdefault@tok@ni\endcsname{\let\PYGdefault@bf=\textbf\def\PYGdefault@tc##1{\textcolor[rgb]{0.60,0.60,0.60}{##1}}}
\expandafter\def\csname PYGdefault@tok@na\endcsname{\def\PYGdefault@tc##1{\textcolor[rgb]{0.49,0.56,0.16}{##1}}}
\expandafter\def\csname PYGdefault@tok@nt\endcsname{\let\PYGdefault@bf=\textbf\def\PYGdefault@tc##1{\textcolor[rgb]{0.00,0.50,0.00}{##1}}}
\expandafter\def\csname PYGdefault@tok@nd\endcsname{\def\PYGdefault@tc##1{\textcolor[rgb]{0.67,0.13,1.00}{##1}}}
\expandafter\def\csname PYGdefault@tok@s\endcsname{\def\PYGdefault@tc##1{\textcolor[rgb]{0.73,0.13,0.13}{##1}}}
\expandafter\def\csname PYGdefault@tok@sd\endcsname{\let\PYGdefault@it=\textit\def\PYGdefault@tc##1{\textcolor[rgb]{0.73,0.13,0.13}{##1}}}
\expandafter\def\csname PYGdefault@tok@si\endcsname{\let\PYGdefault@bf=\textbf\def\PYGdefault@tc##1{\textcolor[rgb]{0.73,0.40,0.53}{##1}}}
\expandafter\def\csname PYGdefault@tok@se\endcsname{\let\PYGdefault@bf=\textbf\def\PYGdefault@tc##1{\textcolor[rgb]{0.73,0.40,0.13}{##1}}}
\expandafter\def\csname PYGdefault@tok@sr\endcsname{\def\PYGdefault@tc##1{\textcolor[rgb]{0.73,0.40,0.53}{##1}}}
\expandafter\def\csname PYGdefault@tok@ss\endcsname{\def\PYGdefault@tc##1{\textcolor[rgb]{0.10,0.09,0.49}{##1}}}
\expandafter\def\csname PYGdefault@tok@sx\endcsname{\def\PYGdefault@tc##1{\textcolor[rgb]{0.00,0.50,0.00}{##1}}}
\expandafter\def\csname PYGdefault@tok@m\endcsname{\def\PYGdefault@tc##1{\textcolor[rgb]{0.40,0.40,0.40}{##1}}}
\expandafter\def\csname PYGdefault@tok@gh\endcsname{\let\PYGdefault@bf=\textbf\def\PYGdefault@tc##1{\textcolor[rgb]{0.00,0.00,0.50}{##1}}}
\expandafter\def\csname PYGdefault@tok@gu\endcsname{\let\PYGdefault@bf=\textbf\def\PYGdefault@tc##1{\textcolor[rgb]{0.50,0.00,0.50}{##1}}}
\expandafter\def\csname PYGdefault@tok@gd\endcsname{\def\PYGdefault@tc##1{\textcolor[rgb]{0.63,0.00,0.00}{##1}}}
\expandafter\def\csname PYGdefault@tok@gi\endcsname{\def\PYGdefault@tc##1{\textcolor[rgb]{0.00,0.63,0.00}{##1}}}
\expandafter\def\csname PYGdefault@tok@gr\endcsname{\def\PYGdefault@tc##1{\textcolor[rgb]{1.00,0.00,0.00}{##1}}}
\expandafter\def\csname PYGdefault@tok@ge\endcsname{\let\PYGdefault@it=\textit}
\expandafter\def\csname PYGdefault@tok@gs\endcsname{\let\PYGdefault@bf=\textbf}
\expandafter\def\csname PYGdefault@tok@gp\endcsname{\let\PYGdefault@bf=\textbf\def\PYGdefault@tc##1{\textcolor[rgb]{0.00,0.00,0.50}{##1}}}
\expandafter\def\csname PYGdefault@tok@go\endcsname{\def\PYGdefault@tc##1{\textcolor[rgb]{0.53,0.53,0.53}{##1}}}
\expandafter\def\csname PYGdefault@tok@gt\endcsname{\def\PYGdefault@tc##1{\textcolor[rgb]{0.00,0.27,0.87}{##1}}}
\expandafter\def\csname PYGdefault@tok@err\endcsname{\def\PYGdefault@bc##1{\setlength{\fboxsep}{0pt}\fcolorbox[rgb]{1.00,0.00,0.00}{1,1,1}{\strut ##1}}}
\expandafter\def\csname PYGdefault@tok@kc\endcsname{\let\PYGdefault@bf=\textbf\def\PYGdefault@tc##1{\textcolor[rgb]{0.00,0.50,0.00}{##1}}}
\expandafter\def\csname PYGdefault@tok@kd\endcsname{\let\PYGdefault@bf=\textbf\def\PYGdefault@tc##1{\textcolor[rgb]{0.00,0.50,0.00}{##1}}}
\expandafter\def\csname PYGdefault@tok@kn\endcsname{\let\PYGdefault@bf=\textbf\def\PYGdefault@tc##1{\textcolor[rgb]{0.00,0.50,0.00}{##1}}}
\expandafter\def\csname PYGdefault@tok@kr\endcsname{\let\PYGdefault@bf=\textbf\def\PYGdefault@tc##1{\textcolor[rgb]{0.00,0.50,0.00}{##1}}}
\expandafter\def\csname PYGdefault@tok@bp\endcsname{\def\PYGdefault@tc##1{\textcolor[rgb]{0.00,0.50,0.00}{##1}}}
\expandafter\def\csname PYGdefault@tok@fm\endcsname{\def\PYGdefault@tc##1{\textcolor[rgb]{0.00,0.00,1.00}{##1}}}
\expandafter\def\csname PYGdefault@tok@vc\endcsname{\def\PYGdefault@tc##1{\textcolor[rgb]{0.10,0.09,0.49}{##1}}}
\expandafter\def\csname PYGdefault@tok@vg\endcsname{\def\PYGdefault@tc##1{\textcolor[rgb]{0.10,0.09,0.49}{##1}}}
\expandafter\def\csname PYGdefault@tok@vi\endcsname{\def\PYGdefault@tc##1{\textcolor[rgb]{0.10,0.09,0.49}{##1}}}
\expandafter\def\csname PYGdefault@tok@vm\endcsname{\def\PYGdefault@tc##1{\textcolor[rgb]{0.10,0.09,0.49}{##1}}}
\expandafter\def\csname PYGdefault@tok@sa\endcsname{\def\PYGdefault@tc##1{\textcolor[rgb]{0.73,0.13,0.13}{##1}}}
\expandafter\def\csname PYGdefault@tok@sb\endcsname{\def\PYGdefault@tc##1{\textcolor[rgb]{0.73,0.13,0.13}{##1}}}
\expandafter\def\csname PYGdefault@tok@sc\endcsname{\def\PYGdefault@tc##1{\textcolor[rgb]{0.73,0.13,0.13}{##1}}}
\expandafter\def\csname PYGdefault@tok@dl\endcsname{\def\PYGdefault@tc##1{\textcolor[rgb]{0.73,0.13,0.13}{##1}}}
\expandafter\def\csname PYGdefault@tok@s2\endcsname{\def\PYGdefault@tc##1{\textcolor[rgb]{0.73,0.13,0.13}{##1}}}
\expandafter\def\csname PYGdefault@tok@sh\endcsname{\def\PYGdefault@tc##1{\textcolor[rgb]{0.73,0.13,0.13}{##1}}}
\expandafter\def\csname PYGdefault@tok@s1\endcsname{\def\PYGdefault@tc##1{\textcolor[rgb]{0.73,0.13,0.13}{##1}}}
\expandafter\def\csname PYGdefault@tok@mb\endcsname{\def\PYGdefault@tc##1{\textcolor[rgb]{0.40,0.40,0.40}{##1}}}
\expandafter\def\csname PYGdefault@tok@mf\endcsname{\def\PYGdefault@tc##1{\textcolor[rgb]{0.40,0.40,0.40}{##1}}}
\expandafter\def\csname PYGdefault@tok@mh\endcsname{\def\PYGdefault@tc##1{\textcolor[rgb]{0.40,0.40,0.40}{##1}}}
\expandafter\def\csname PYGdefault@tok@mi\endcsname{\def\PYGdefault@tc##1{\textcolor[rgb]{0.40,0.40,0.40}{##1}}}
\expandafter\def\csname PYGdefault@tok@il\endcsname{\def\PYGdefault@tc##1{\textcolor[rgb]{0.40,0.40,0.40}{##1}}}
\expandafter\def\csname PYGdefault@tok@mo\endcsname{\def\PYGdefault@tc##1{\textcolor[rgb]{0.40,0.40,0.40}{##1}}}
\expandafter\def\csname PYGdefault@tok@ch\endcsname{\let\PYGdefault@it=\textit\def\PYGdefault@tc##1{\textcolor[rgb]{0.25,0.50,0.50}{##1}}}
\expandafter\def\csname PYGdefault@tok@cm\endcsname{\let\PYGdefault@it=\textit\def\PYGdefault@tc##1{\textcolor[rgb]{0.25,0.50,0.50}{##1}}}
\expandafter\def\csname PYGdefault@tok@cpf\endcsname{\let\PYGdefault@it=\textit\def\PYGdefault@tc##1{\textcolor[rgb]{0.25,0.50,0.50}{##1}}}
\expandafter\def\csname PYGdefault@tok@c1\endcsname{\let\PYGdefault@it=\textit\def\PYGdefault@tc##1{\textcolor[rgb]{0.25,0.50,0.50}{##1}}}
\expandafter\def\csname PYGdefault@tok@cs\endcsname{\let\PYGdefault@it=\textit\def\PYGdefault@tc##1{\textcolor[rgb]{0.25,0.50,0.50}{##1}}}

\makeatother

\makeatletter
\def\PYG@reset{\let\PYG@it=\relax \let\PYG@bf=\relax%
    \let\PYG@ul=\relax \let\PYG@tc=\relax%
    \let\PYG@bc=\relax \let\PYG@ff=\relax}
\def\PYG@tok#1{\csname PYG@tok@#1\endcsname}
\def\PYG@toks#1+{\ifx\relax#1\empty\else%
    \PYG@tok{#1}\expandafter\PYG@toks\fi}
\def\PYG@do#1{\PYG@bc{\PYG@tc{\PYG@ul{%
    \PYG@it{\PYG@bf{\PYG@ff{#1}}}}}}}
\def\PYG#1#2{\PYG@reset\PYG@toks#1+\relax+\PYG@do{#2}}

\expandafter\def\csname PYG@tok@w\endcsname{\def\PYG@tc##1{\textcolor[rgb]{0.73,0.73,0.73}{##1}}}
\expandafter\def\csname PYG@tok@c\endcsname{\let\PYG@it=\textit\def\PYG@tc##1{\textcolor[rgb]{0.25,0.50,0.50}{##1}}}
\expandafter\def\csname PYG@tok@cp\endcsname{\def\PYG@tc##1{\textcolor[rgb]{0.74,0.48,0.00}{##1}}}
\expandafter\def\csname PYG@tok@k\endcsname{\let\PYG@bf=\textbf\def\PYG@tc##1{\textcolor[rgb]{0.00,0.50,0.00}{##1}}}
\expandafter\def\csname PYG@tok@kp\endcsname{\def\PYG@tc##1{\textcolor[rgb]{0.00,0.50,0.00}{##1}}}
\expandafter\def\csname PYG@tok@kt\endcsname{\def\PYG@tc##1{\textcolor[rgb]{0.69,0.00,0.25}{##1}}}
\expandafter\def\csname PYG@tok@o\endcsname{\def\PYG@tc##1{\textcolor[rgb]{0.40,0.40,0.40}{##1}}}
\expandafter\def\csname PYG@tok@ow\endcsname{\let\PYG@bf=\textbf\def\PYG@tc##1{\textcolor[rgb]{0.67,0.13,1.00}{##1}}}
\expandafter\def\csname PYG@tok@nb\endcsname{\def\PYG@tc##1{\textcolor[rgb]{0.00,0.50,0.00}{##1}}}
\expandafter\def\csname PYG@tok@nf\endcsname{\def\PYG@tc##1{\textcolor[rgb]{0.00,0.00,1.00}{##1}}}
\expandafter\def\csname PYG@tok@nc\endcsname{\let\PYG@bf=\textbf\def\PYG@tc##1{\textcolor[rgb]{0.00,0.00,1.00}{##1}}}
\expandafter\def\csname PYG@tok@nn\endcsname{\let\PYG@bf=\textbf\def\PYG@tc##1{\textcolor[rgb]{0.00,0.00,1.00}{##1}}}
\expandafter\def\csname PYG@tok@ne\endcsname{\let\PYG@bf=\textbf\def\PYG@tc##1{\textcolor[rgb]{0.82,0.25,0.23}{##1}}}
\expandafter\def\csname PYG@tok@nv\endcsname{\def\PYG@tc##1{\textcolor[rgb]{0.10,0.09,0.49}{##1}}}
\expandafter\def\csname PYG@tok@no\endcsname{\def\PYG@tc##1{\textcolor[rgb]{0.53,0.00,0.00}{##1}}}
\expandafter\def\csname PYG@tok@nl\endcsname{\def\PYG@tc##1{\textcolor[rgb]{0.63,0.63,0.00}{##1}}}
\expandafter\def\csname PYG@tok@ni\endcsname{\let\PYG@bf=\textbf\def\PYG@tc##1{\textcolor[rgb]{0.60,0.60,0.60}{##1}}}
\expandafter\def\csname PYG@tok@na\endcsname{\def\PYG@tc##1{\textcolor[rgb]{0.49,0.56,0.16}{##1}}}
\expandafter\def\csname PYG@tok@nt\endcsname{\let\PYG@bf=\textbf\def\PYG@tc##1{\textcolor[rgb]{0.00,0.50,0.00}{##1}}}
\expandafter\def\csname PYG@tok@nd\endcsname{\def\PYG@tc##1{\textcolor[rgb]{0.67,0.13,1.00}{##1}}}
\expandafter\def\csname PYG@tok@s\endcsname{\def\PYG@tc##1{\textcolor[rgb]{0.73,0.13,0.13}{##1}}}
\expandafter\def\csname PYG@tok@sd\endcsname{\let\PYG@it=\textit\def\PYG@tc##1{\textcolor[rgb]{0.73,0.13,0.13}{##1}}}
\expandafter\def\csname PYG@tok@si\endcsname{\let\PYG@bf=\textbf\def\PYG@tc##1{\textcolor[rgb]{0.73,0.40,0.53}{##1}}}
\expandafter\def\csname PYG@tok@se\endcsname{\let\PYG@bf=\textbf\def\PYG@tc##1{\textcolor[rgb]{0.73,0.40,0.13}{##1}}}
\expandafter\def\csname PYG@tok@sr\endcsname{\def\PYG@tc##1{\textcolor[rgb]{0.73,0.40,0.53}{##1}}}
\expandafter\def\csname PYG@tok@ss\endcsname{\def\PYG@tc##1{\textcolor[rgb]{0.10,0.09,0.49}{##1}}}
\expandafter\def\csname PYG@tok@sx\endcsname{\def\PYG@tc##1{\textcolor[rgb]{0.00,0.50,0.00}{##1}}}
\expandafter\def\csname PYG@tok@m\endcsname{\def\PYG@tc##1{\textcolor[rgb]{0.40,0.40,0.40}{##1}}}
\expandafter\def\csname PYG@tok@gh\endcsname{\let\PYG@bf=\textbf\def\PYG@tc##1{\textcolor[rgb]{0.00,0.00,0.50}{##1}}}
\expandafter\def\csname PYG@tok@gu\endcsname{\let\PYG@bf=\textbf\def\PYG@tc##1{\textcolor[rgb]{0.50,0.00,0.50}{##1}}}
\expandafter\def\csname PYG@tok@gd\endcsname{\def\PYG@tc##1{\textcolor[rgb]{0.63,0.00,0.00}{##1}}}
\expandafter\def\csname PYG@tok@gi\endcsname{\def\PYG@tc##1{\textcolor[rgb]{0.00,0.63,0.00}{##1}}}
\expandafter\def\csname PYG@tok@gr\endcsname{\def\PYG@tc##1{\textcolor[rgb]{1.00,0.00,0.00}{##1}}}
\expandafter\def\csname PYG@tok@ge\endcsname{\let\PYG@it=\textit}
\expandafter\def\csname PYG@tok@gs\endcsname{\let\PYG@bf=\textbf}
\expandafter\def\csname PYG@tok@gp\endcsname{\let\PYG@bf=\textbf\def\PYG@tc##1{\textcolor[rgb]{0.00,0.00,0.50}{##1}}}
\expandafter\def\csname PYG@tok@go\endcsname{\def\PYG@tc##1{\textcolor[rgb]{0.53,0.53,0.53}{##1}}}
\expandafter\def\csname PYG@tok@gt\endcsname{\def\PYG@tc##1{\textcolor[rgb]{0.00,0.27,0.87}{##1}}}
\expandafter\def\csname PYG@tok@err\endcsname{\def\PYG@bc##1{\setlength{\fboxsep}{0pt}\fcolorbox[rgb]{1.00,0.00,0.00}{1,1,1}{\strut ##1}}}
\expandafter\def\csname PYG@tok@kc\endcsname{\let\PYG@bf=\textbf\def\PYG@tc##1{\textcolor[rgb]{0.00,0.50,0.00}{##1}}}
\expandafter\def\csname PYG@tok@kd\endcsname{\let\PYG@bf=\textbf\def\PYG@tc##1{\textcolor[rgb]{0.00,0.50,0.00}{##1}}}
\expandafter\def\csname PYG@tok@kn\endcsname{\let\PYG@bf=\textbf\def\PYG@tc##1{\textcolor[rgb]{0.00,0.50,0.00}{##1}}}
\expandafter\def\csname PYG@tok@kr\endcsname{\let\PYG@bf=\textbf\def\PYG@tc##1{\textcolor[rgb]{0.00,0.50,0.00}{##1}}}
\expandafter\def\csname PYG@tok@bp\endcsname{\def\PYG@tc##1{\textcolor[rgb]{0.00,0.50,0.00}{##1}}}
\expandafter\def\csname PYG@tok@fm\endcsname{\def\PYG@tc##1{\textcolor[rgb]{0.00,0.00,1.00}{##1}}}
\expandafter\def\csname PYG@tok@vc\endcsname{\def\PYG@tc##1{\textcolor[rgb]{0.10,0.09,0.49}{##1}}}
\expandafter\def\csname PYG@tok@vg\endcsname{\def\PYG@tc##1{\textcolor[rgb]{0.10,0.09,0.49}{##1}}}
\expandafter\def\csname PYG@tok@vi\endcsname{\def\PYG@tc##1{\textcolor[rgb]{0.10,0.09,0.49}{##1}}}
\expandafter\def\csname PYG@tok@vm\endcsname{\def\PYG@tc##1{\textcolor[rgb]{0.10,0.09,0.49}{##1}}}
\expandafter\def\csname PYG@tok@sa\endcsname{\def\PYG@tc##1{\textcolor[rgb]{0.73,0.13,0.13}{##1}}}
\expandafter\def\csname PYG@tok@sb\endcsname{\def\PYG@tc##1{\textcolor[rgb]{0.73,0.13,0.13}{##1}}}
\expandafter\def\csname PYG@tok@sc\endcsname{\def\PYG@tc##1{\textcolor[rgb]{0.73,0.13,0.13}{##1}}}
\expandafter\def\csname PYG@tok@dl\endcsname{\def\PYG@tc##1{\textcolor[rgb]{0.73,0.13,0.13}{##1}}}
\expandafter\def\csname PYG@tok@s2\endcsname{\def\PYG@tc##1{\textcolor[rgb]{0.73,0.13,0.13}{##1}}}
\expandafter\def\csname PYG@tok@sh\endcsname{\def\PYG@tc##1{\textcolor[rgb]{0.73,0.13,0.13}{##1}}}
\expandafter\def\csname PYG@tok@s1\endcsname{\def\PYG@tc##1{\textcolor[rgb]{0.73,0.13,0.13}{##1}}}
\expandafter\def\csname PYG@tok@mb\endcsname{\def\PYG@tc##1{\textcolor[rgb]{0.40,0.40,0.40}{##1}}}
\expandafter\def\csname PYG@tok@mf\endcsname{\def\PYG@tc##1{\textcolor[rgb]{0.40,0.40,0.40}{##1}}}
\expandafter\def\csname PYG@tok@mh\endcsname{\def\PYG@tc##1{\textcolor[rgb]{0.40,0.40,0.40}{##1}}}
\expandafter\def\csname PYG@tok@mi\endcsname{\def\PYG@tc##1{\textcolor[rgb]{0.40,0.40,0.40}{##1}}}
\expandafter\def\csname PYG@tok@il\endcsname{\def\PYG@tc##1{\textcolor[rgb]{0.40,0.40,0.40}{##1}}}
\expandafter\def\csname PYG@tok@mo\endcsname{\def\PYG@tc##1{\textcolor[rgb]{0.40,0.40,0.40}{##1}}}
\expandafter\def\csname PYG@tok@ch\endcsname{\let\PYG@it=\textit\def\PYG@tc##1{\textcolor[rgb]{0.25,0.50,0.50}{##1}}}
\expandafter\def\csname PYG@tok@cm\endcsname{\let\PYG@it=\textit\def\PYG@tc##1{\textcolor[rgb]{0.25,0.50,0.50}{##1}}}
\expandafter\def\csname PYG@tok@cpf\endcsname{\let\PYG@it=\textit\def\PYG@tc##1{\textcolor[rgb]{0.25,0.50,0.50}{##1}}}
\expandafter\def\csname PYG@tok@c1\endcsname{\let\PYG@it=\textit\def\PYG@tc##1{\textcolor[rgb]{0.25,0.50,0.50}{##1}}}
\expandafter\def\csname PYG@tok@cs\endcsname{\let\PYG@it=\textit\def\PYG@tc##1{\textcolor[rgb]{0.25,0.50,0.50}{##1}}}

\def\PYGZhy{\char`\-}
\def\PYGZsq{\char`\'}

\makeatother

\usepackage{algorithmic}
\ifpdf
  \DeclareGraphicsExtensions{.eps,.pdf,.png,.jpg}
\else
  \DeclareGraphicsExtensions{.eps}
\fi

\usepackage{tikz}

\numberwithin{theorem}{section}

\newcommand{\TheTitle}{Chaos Game Representation} 
\newcommand{\TheAuthors}{E.~Y.~S.~Chan, and R.~M.~Corless}

\headers{\TheTitle}{\TheAuthors}

\title{{\TheTitle}\thanks{Submitted to the editors DATE.}}

\author{
  Eunice Y. S. Chan\thanks{Department of Anesthesia and Perioperative Medicine, MEDICI Centre, Western University, London, Ontario, Canada
    (\email{echan295@uwo.ca}).}
  \and
  Robert M. Corless\thanks{Ontario Research Center for Computer Algebra, School of Mathematical and Statistical Sciences, and Rotman Institute of Philosophy, Western University, London, Ontario, Canada (\email{rcorless@uwo.ca}).}
}

\usepackage{amsopn}

\ifpdf
\hypersetup{
  pdftitle={\TheTitle},
  pdfauthor={\TheAuthors}
}
\fi

\usepackage{tikz}

\begin{document}

\maketitle

\begin{abstract}
	The chaos game representation (CGR) is an interesting method to visualize one-dimensional sequences. In this paper, we show how to construct a chaos game representation. The applications mentioned here are biological, in which CGR was able to uncover patterns in DNA or proteins that were previously unknown. We also show how CGR might be introduced in the classroom, either in a modelling course or in a dynamical systems course. Some sequences that are tested are taken from the Online Encyclopedia of Integer Sequences, and others are taken from sequences that arose mainly from a course in experimental mathematics.
\end{abstract}

\begin{keywords}
	Chaos game representation, iterated function systems, partial quotients of continued fractions, sequences
\end{keywords}

\begin{AMS}
  
\end{AMS}

\section{Introduction}

Finding hidden patterns in long sequences can be both difficult and valuable. Representing these sequence in a visual way can often help. The so-called chaos game representation (CGR) of a sequence of integers is a particularly useful technique, that visualizes a one-dimensional sequence in a two-dimensional space. The CGR is presented as a scatter plot (most frequently square), in which each corner represents an element that appears in the sequence. The results of these CGRs can look very fractal-like, but even so can be visually recognizable and distinguishable. The distinguishing characteristics can be made quantitative, with a good notion of ``distance between images.''

Many applications, such as analysis of DNA sequences~\cite{goldman1993nucleotide, jeffrey1990chaos, jeffrey1992chaos, karamichalis2015investigation} and protein structure~\cite{basu1997chaos, fiser1994chaos}, have shown the usefulness of CGR; we will discuss these applications briefly in section~\ref{sec:DNA}. Before that, in section~\ref{sec:ifs}, we will look at a random iteration algorithm for creating fractals. In section~\ref{sec:math}, we will apply CGR to simple abstract mathematical sequences, including the digits of $\pi$ and the partial quotients of continued fractions. Additionally, we will look at certain ``distances between images" that are produced. Finally, we will explain how some of the patterns arise, and what these depictions can tell us about the sequences.

One important pedagogical purpose of this module is to provoke a discussion about randomness, or what it means to be random. In this paper, we at first use the word ``random" very loosely, on purpose.

\begin{displayquote}
	``What is probability? I asked myself this question many years ago, and found that various authors gave different answers. I found that there were several main schools of thought with many variations." \vbox{\hfil--- Richard W.~Hamming\footnote{Bio: \url{http://mathshistory.st-andrews.ac.uk/Biographies/Hamming.html}} \cite{hamming1991art}}
\end{displayquote}

\begin{figure}
	\centering
	\includegraphics[width=0.5\textwidth]{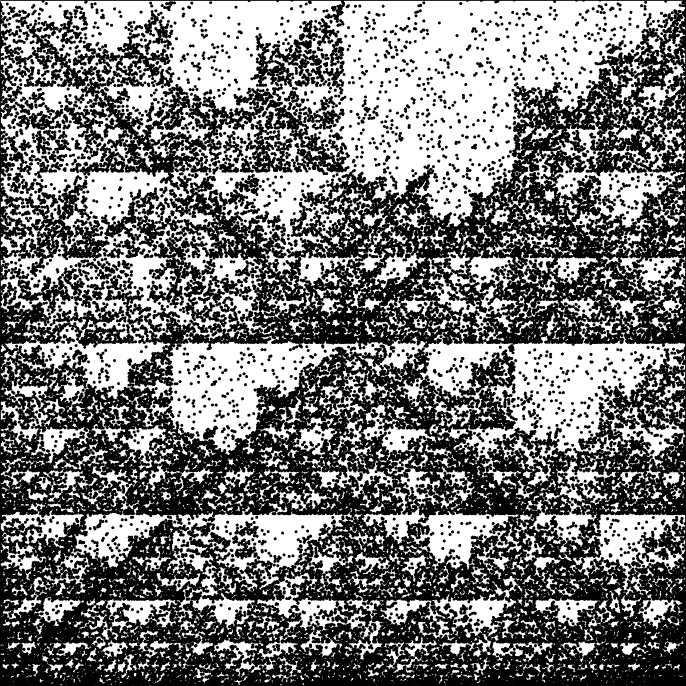}
	\caption{DNA of human beta globin region on chromosome 11 (HUMHBB)---73308 bps}
	\label{fig:cgr_globin}
\end{figure}


\section{Creating Fractals using Iterated Function Systems (IFS)} \label{sec:ifs}
The random iteration algorithm for creating pictures of fractals uses the fixed attracting set of an \emph{iterated function system} (IFS). Before defining IFS, we will first define some necessary preliminaries.

\subsection{Affine Transformation in the Euclidean Plane}
\begin{figure}
	\centering
	\begin{subfigure}[b]{0.45\textwidth}
		\includegraphics[width=\textwidth]{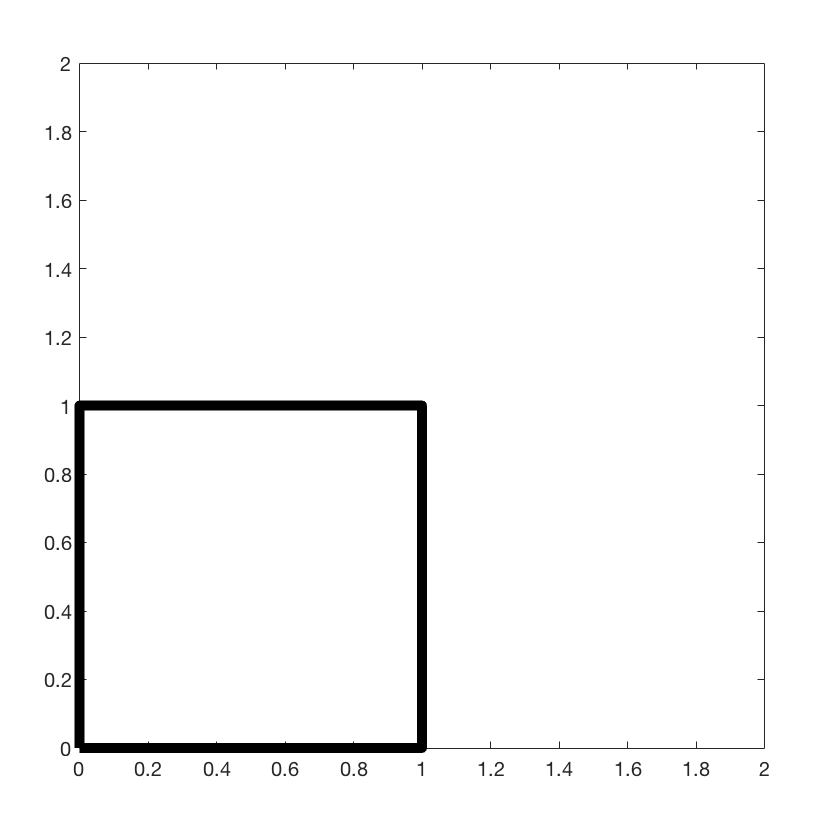}
		\caption{Shape before transformation}
	\end{subfigure}
	\quad
	\begin{subfigure}[b]{0.45\textwidth}
		\includegraphics[width=\textwidth]{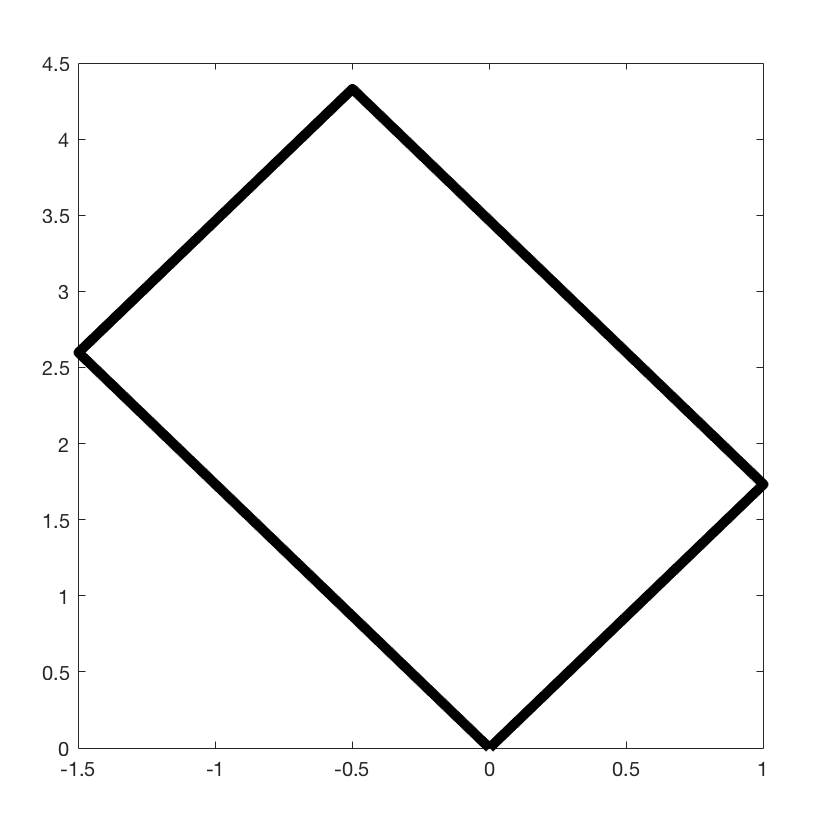}
		\caption{Shape after transformation}
	\end{subfigure}
	\caption{Transformation of a unit square using the transformation matrix $A$ (as shown in equation \eqref{eqn:tmatrix}), where $r_{1} = 2$, $r_{2} = 3$, $\theta_{1} = \sfrac{\pi}{3}$, and $\theta_{2} = \sfrac{\pi}{6}$.}
	\label{fig:transform}
\end{figure}

One can write a two-dimensional \emph{affine} transformation in the Euclidean plane $w: \mathbb{R}^{2} \to \mathbb{R}^{2}$ in the form
\begin{equation}
	w(x_{1}, x_{2}) = \left(ax_{1} + bx_{2} + e, cx_{1} + dx_{2} + f \right)
\end{equation}
where $a$, $b$, $c$, $d$, $e$, and $f$ are real numbers~\cite{barnsley2014fractals}. We will use the following more compact notation:
\begin{align}
	w(x) &= 
	w
	\begin{pmatrix}
		x_{1} \\
		x_{2}
	\end{pmatrix}
	=
	\begin{bmatrix}
		a & b \\
		c & d
	\end{bmatrix}
	\begin{bmatrix}
		x_{1} \\
		x_{2}
	\end{bmatrix}
	+
	\begin{bmatrix}
		e \\
		f
	\end{bmatrix} \\ \label{eqn:ifs} 
	&= Ax + t \nonumber \>.
\end{align}
The matrix $A$ is a $2 \times 2$ real matrix and can always be written in the form
\begin{equation}
	\begin{bmatrix}
		a & b \\
		c & d
	\end{bmatrix}
	=
	\begin{bmatrix}
		r_{1}\cos\theta_{1} & -r_{2}\sin\theta_{2} \\
		r_{1}\sin\theta_{1} & \phantom{-}r_{2}\cos\theta_{2}
	\end{bmatrix} \>,
	\label{eqn:tmatrix}
\end{equation}
where $\left(r_{1}, \theta_{1}\right)$ are the polar coordinates of the point $\left(a, c\right)$ and $\left(r_{2}, \left(\theta_{2} + \sfrac{\pi}{2}\right)\right)$ are the polar coordinates of the point $\left(b, d\right)$. An example of a \emph{linear} transformation
\begin{equation}
	\begin{bmatrix}
		x_{1} \\
		x_{2}
	\end{bmatrix}
	\to
	A
	\begin{bmatrix}
		x_{1} \\
		x_{2}
	\end{bmatrix}
\end{equation}
in $\mathbb{R}^{2}$ is shown in figure \ref{fig:transform}. We can see from this figure that the original shape remains as a parallelogram, but now has changed in size and in rotation. This is what we mean by \emph{linear} transformation.

\subsection{Iterated Function System}
An \emph{iterated function system} (IFS)\footnote{To be more specific, this is a \emph{hyperbolic} iterated function system. However, in practice, the word ``hyperbolic" is sometimes dropped in practice~\cite[p.~82]{barnsley2014fractals}.} is a finite set of contraction mappings on a metric space.  These can be used to create pictures of fractals. These pictures essentially show what is called the attractor of the IFS. We will be defining what an attractor is in Section~\ref{sec:chaos}. In this paper, we will only be looking at the form where each contraction mapping is an affine transformation on the metric space $\mathbb{R}^2$. 

We will illustrate the algorithms for an IFS as an example. Consider the following maps:
\begin{align}
	w_1(x, y) &=
	\left[
		\begin{array}{cc}
			\sfrac{1}{2} & 0 \\
			0 & \sfrac{1}{2}
		\end{array}
	\right]
	\left[
		\begin{array}{c}
			x \\
			y
		\end{array}
	\right] 
	+
	\left[
		\begin{array}{c}
			0 \\
			0
		\end{array}
	\right] \>, \nonumber \\
	w_2(x, y) &=
	\left[
		\begin{array}{cc}
			\sfrac{1}{2} & 0 \\
			0 & \sfrac{1}{2}
		\end{array}
	\right]
	\left[
		\begin{array}{c}
			x \\
			y
		\end{array}
	\right]
	+
	\left[
		\begin{array}{c}
			0 \\
			\sfrac{1}{2}
		\end{array}
	\right] \>, \\
	w_3(x, y) &=
	\left[
		\begin{array}{cc}
			\sfrac{1}{2} & 0 \\
			0 & \sfrac{1}{2}
		\end{array}
	\right]
	\left[
		\begin{array}{c}
			x \\
			y
		\end{array}
	\right]
	+
	\left[
		\begin{array}{c}
			\sfrac{1}{2} \\
			\sfrac{1}{2}
		\end{array}
	\right] \nonumber \>,
	\label{eqn:sierpinski}
\end{align}
each with the \emph{probability factor} of $\sfrac{1}{3}$. These values can also be displayed as a table, where the order of the coefficients $a$ through $f$ corresponds to the order presented in equation \eqref{eqn:ifs} and $p$ represents the probability factor. The corresponding table for equation \eqref{eqn:sierpinski} can be found in table \ref{tab:IFS_sierpinski}.

\begin{table}[ht]
    \centering
    \caption{Table representing the IFS from equation \eqref{eqn:sierpinski}}
	\begin{tabular}{c|ccccccc}
		$w$ & $a$ & $b$ & $c$ & $d$ & $e$ & $f$ & $p$\\
		\hline
		$1$ & $\sfrac{1}{2}$ & $0$ & $0$ & $\sfrac{1}{2}$ & $0$ & $0$ & $\sfrac{1}{3}$\\
		$2$ & $\sfrac{1}{2}$ & $0$ & $0$ & $\sfrac{1}{2}$ & $0$ & $\sfrac{1}{2}$ & $\sfrac{1}{3}$\\
		$3$ & $\sfrac{1}{2}$ & $0$ & $0$ & $\sfrac{1}{2}$ & $\sfrac{1}{2}$ & $\sfrac{1}{2}$ & $\sfrac{1}{3}$\\
	\end{tabular}
	\label{tab:IFS_sierpinski}
\end{table}

\begin{figure}[ht]
    \centering
    \includegraphics[width= 0.5\textwidth]{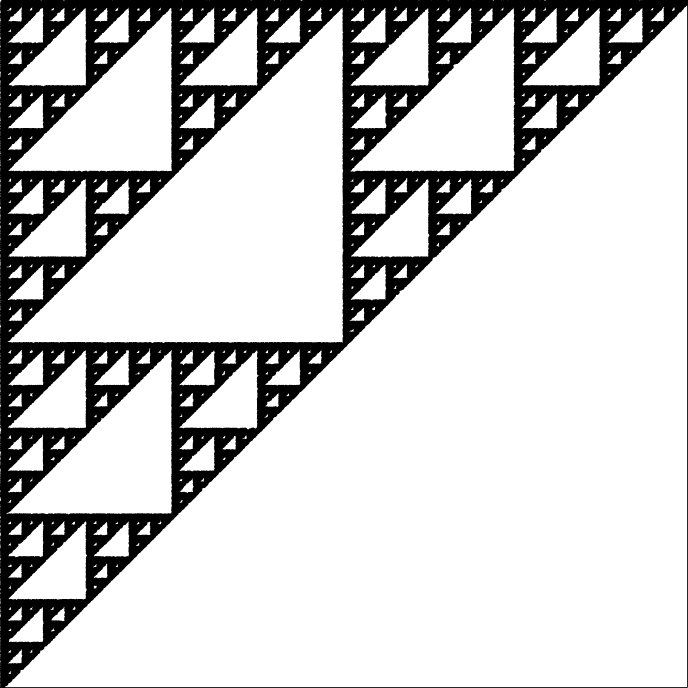}
    \caption{Results from Listing \ref{alg:sierpinski}. We see visible regularity arising from a program that uses randomness.}
    \label{fig:ifs_Sier}
\end{figure}

To create a fractal picture using IFS, we first choose a starting point $(x_0, y_0)$. We then \emph{randomly} choose a map from our IFS and evaluate it at our starting point $(x_0, y_0)$ to get our next point $(x_1, y_1)$. We do this again many times (determined by the user) until a pattern (usually fractal) appears. The \textsc{Matlab} code shown in Listing \ref{alg:sierpinski} computes and plots 50000 points corresponding to the IFS from table \ref{tab:IFS_sierpinski}. To choose which map to use for each iteration, we used a uniform random number generator \texttt{randi} (built-in function in \textsc{Matlab}) to choose a number between (and including) 0 and 2. Each value corresponds to a map in our IFS. The picture produced from Listing \ref{alg:sierpinski} is shown in Figure \ref{fig:ifs_Sier}. It is clear that this figure shows a fractal, which appears to be the Sierpinski triangle~\cite{barnsley2014fractals, feldman2012chaos}, where its three vertices are located at $(0, 0)$, $(0, 1)$ and $(1, 1)$.

\begin{listing}[ht]
\caption{\textsc{Matlab} Code for an IFS generating a Sierpinski triangle}
\label{alg:sierpinski}
\begin{Verbatim}[commandchars=\\\{\}]
\PYG{n}{n} \PYG{p}{=} \PYG{l+m+mi}{50000}\PYG{p}{;}
\PYG{n}{pts} \PYG{p}{=} \PYG{n+nb}{zeros}\PYG{p}{(}\PYG{l+m+mi}{2}\PYG{p}{,} \PYG{n}{n}\PYG{p}{);}
\PYG{n}{rseq} \PYG{p}{=} \PYG{n}{randi}\PYG{p}{([}\PYG{l+m+mi}{0}\PYG{p}{,} \PYG{l+m+mi}{2}\PYG{p}{],} \PYG{l+m+mi}{1}\PYG{p}{,} \PYG{n}{n}\PYG{o}{\PYGZhy{}}\PYG{l+m+mi}{1}\PYG{p}{);}
\PYG{k}{for} \PYG{n+nb}{i} \PYG{p}{=} \PYG{l+m+mi}{2}\PYG{p}{:}\PYG{n}{n}
	\PYG{k}{if} \PYG{n}{rseq}\PYG{p}{(}\PYG{n+nb}{i}\PYG{o}{\PYGZhy{}}\PYG{l+m+mi}{1}\PYG{p}{)} \PYG{o}{==} \PYG{l+m+mi}{0}
		\PYG{n}{pts}\PYG{p}{(:,} \PYG{n+nb}{i}\PYG{p}{)} \PYG{p}{=} \PYG{p}{[}\PYG{l+m+mf}{0.5}\PYG{p}{,} \PYG{l+m+mi}{0}\PYG{p}{;} \PYG{l+m+mi}{0}\PYG{p}{,} \PYG{l+m+mf}{0.5}\PYG{p}{]}\PYG{o}{*}\PYG{n}{pts}\PYG{p}{(:,} \PYG{n+nb}{i}\PYG{o}{\PYGZhy{}}\PYG{l+m+mi}{1}\PYG{p}{)} \PYG{o}{+} \PYG{p}{[}\PYG{l+m+mi}{0}\PYG{p}{;} \PYG{l+m+mi}{0}\PYG{p}{];}
	\PYG{k}{elseif} \PYG{n}{rseq}\PYG{p}{(}\PYG{n+nb}{i}\PYG{o}{\PYGZhy{}}\PYG{l+m+mi}{1}\PYG{p}{)} \PYG{o}{==} \PYG{l+m+mi}{1}
		\PYG{n}{pts}\PYG{p}{(:,} \PYG{n+nb}{i}\PYG{p}{)} \PYG{p}{=} \PYG{p}{[}\PYG{l+m+mf}{0.5}\PYG{p}{,} \PYG{l+m+mi}{0}\PYG{p}{;} \PYG{l+m+mi}{0}\PYG{p}{,} \PYG{l+m+mf}{0.5}\PYG{p}{]}\PYG{o}{*}\PYG{n}{pts}\PYG{p}{(:,} \PYG{n+nb}{i}\PYG{o}{\PYGZhy{}}\PYG{l+m+mi}{1}\PYG{p}{)} \PYG{o}{+} \PYG{p}{[}\PYG{l+m+mi}{0}\PYG{p}{;} \PYG{l+m+mf}{0.5}\PYG{p}{];}
	\PYG{k}{else}
		\PYG{n}{pts}\PYG{p}{(:,} \PYG{n+nb}{i}\PYG{p}{)} \PYG{p}{=} \PYG{p}{[}\PYG{l+m+mf}{0.5}\PYG{p}{,} \PYG{l+m+mi}{0}\PYG{p}{;} \PYG{l+m+mi}{0}\PYG{p}{,} \PYG{l+m+mf}{0.5}\PYG{p}{]}\PYG{o}{*}\PYG{n}{pts}\PYG{p}{(:,} \PYG{n+nb}{i}\PYG{o}{\PYGZhy{}}\PYG{l+m+mi}{1}\PYG{p}{)} \PYG{o}{+} \PYG{p}{[}\PYG{l+m+mf}{0.5}\PYG{p}{;} \PYG{l+m+mf}{0.5}\PYG{p}{];}
	\PYG{k}{end}
\PYG{k}{end}
\PYG{n}{plot}\PYG{p}{(}\PYG{n}{pts}\PYG{p}{(}\PYG{l+m+mi}{1}\PYG{p}{,} \PYG{p}{:),} \PYG{n}{pts}\PYG{p}{(}\PYG{l+m+mi}{2}\PYG{p}{,} \PYG{p}{:),} \PYG{l+s}{\PYGZsq{}k.\PYGZsq{}}\PYG{p}{,} \PYG{l+s}{\PYGZsq{}MarkerSize\PYGZsq{}}\PYG{p}{,} \PYG{l+m+mi}{1}\PYG{p}{);}
\PYG{n}{axis}\PYG{p}{(}\PYG{l+s}{\PYGZsq{}square\PYGZsq{}}\PYG{p}{);}
\PYG{n}{set}\PYG{p}{(}\PYG{n}{gca}\PYG{p}{,} \PYG{l+s}{\PYGZsq{}xtick\PYGZsq{}}\PYG{p}{,} \PYG{p}{[]);}
\PYG{n}{set}\PYG{p}{(}\PYG{n}{gca}\PYG{p}{,} \PYG{l+s}{\PYGZsq{}ytick\PYGZsq{}}\PYG{p}{,} \PYG{p}{[]);}
\end{Verbatim}
\end{listing}

\begin{figure}[th]
	\centering
	\begin{subfigure}[b]{0.475\textwidth}
		\includegraphics[width=\textwidth]{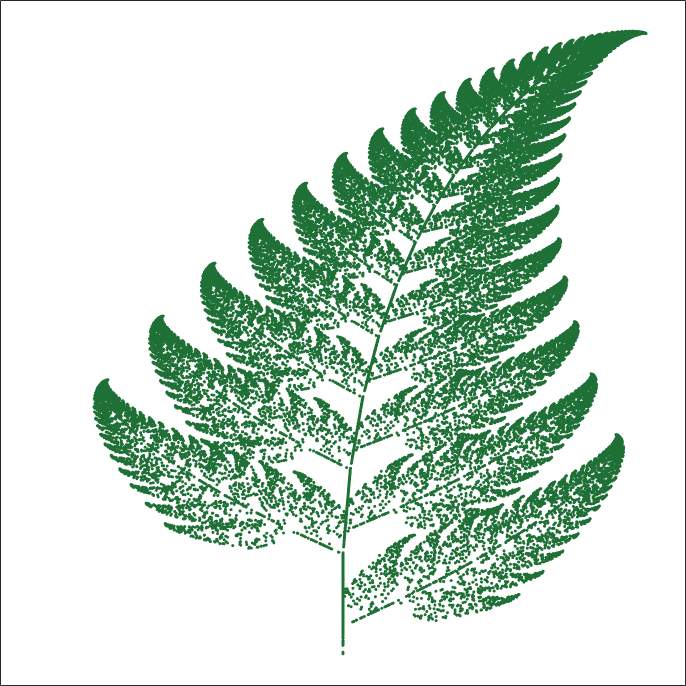}
		\caption{Barnsley's Fern (n = 50000)}
		\label{fig:barnsley_fern}
	\end{subfigure}
	\quad
	\begin{subfigure}[b]{0.475\textwidth}
		\includegraphics[width=\textwidth]{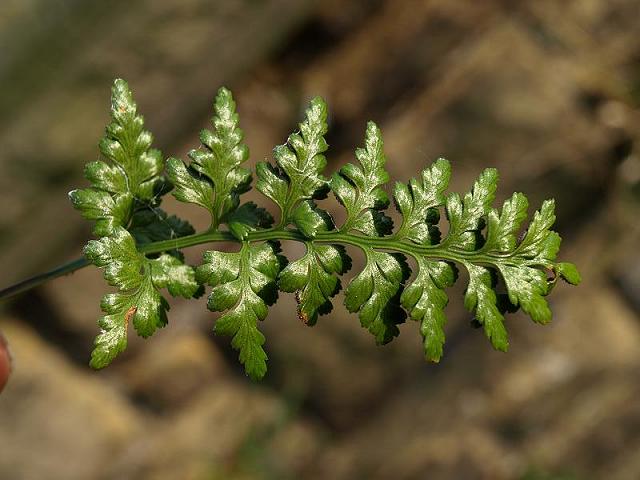}
		\caption{Black spleenwort\footnote{\url{http://www.aphotoflora.com/images/pteridophyta_aspleniaceae_ferns/asplenium_adiantum-nigrum_leaf_27-09-08.jpg}}}
		\label{fig:black_spleenwort}
	\end{subfigure}
	\caption{Comparison of Barnsley's Fern to the black spleenwort, which is what the IFS was based on.}
	\label{fig:fern}
\end{figure}

Iterated function systems can be quite versatile; a set of several attractors can be use to create different shapes and even can resemble real-life objects. An example of this is Barnsley's fern, which resembles the Black Spleenwort~\cite{barnsley2014fractals} (see Figure \ref{fig:black_spleenwort}). The map is displayed in the following table:

\begin{center}
	\begin{tabular}{c|ccccccc}
		$w$ & $a$ & $b$ & $c$ & $d$ & $e$ & $f$ & $p$\\
		\hline
		$1$ & $\phantom{-}0.00$ & $\phantom{-}0.00$ & $\phantom{-}0.00$ & $0.16$ & $0.00$ & $0.00$ & $0.01$\\
		$2$ & $\phantom{-}0.85$ & $\phantom{-}0.04$ & $-0.04$ & $0.85$ & $0.00$ & $1.60$ & $0.85$\\
		$3$ & $\phantom{-}0.20$ & $-0.26$ & $\phantom{-}0.23$ & $0.22$ & $0.00$ & $1.60$ & $0.07$\\
		$4$ & $-0.15$ & $\phantom{-}0.28$ & $\phantom{-}0.26$ & $0.24$ & $0.00$ & $0.44$ & $0.07$
	\end{tabular}
\end{center}
The plot produced by this IFS is shown in figure \ref{fig:barnsley_fern}.

\section{Chaos game} \label{sec:chaos}
CGR is based on a technique from chaotic dynamics called the ``chaos game," popularized by Barnsley~\cite{barnsley2014fractals} in 1988. The chaos game is an algorithm which allows one to produce pictures of fractal structures~\cite{jeffrey1990chaos}. In its simplest form, this can be demonstrated with a piece of paper and pencil using the following steps (adapted from \cite{jeffrey1990chaos}):
\begin{enumerate}
	\item On a piece of paper, draw three points spaced out across the page; they do not necessarily need to be spaced out evenly to form an equilateral triangle though they might be.
	\item Label one of the points with ``1, 2,'' the other point with ``3, 4,'' and the last point with ``5, 6.''
	\item Draw a point anywhere on the page: this will be the starting point.
	\item Roll a six-sided die. Draw an imaginary line from the current dot towards to the point with the corresponding number, and put a dot half way on the line. This becomes the new ``current dot.''
	\item Repeat step 4 (until you get bored).
\end{enumerate}
Obviously, one would not want to plot a large number of points by hand\footnote{We tried this, though not in class. We got up to 80 points on one figure, using the die roller on \url{random.org}. It's kind of fun, in a physical ``drawing and measuring" way, but it's hard not to make mistakes. It could make a useful ``active learning" exercise and we might try it in class in the future.}; instead, this can be done by computer; the code is shown in Listing \ref{alg:triangle}. For this algorithm, we assigned $(0, 0)$, $(0, 1)$ and $(1, 1)$ as the vertices of the triangle; they each are labelled with the value $0$, $1$, and $2$, respectively. Instead of using a die to decide which vertex to choose, similarly to Listing \ref{alg:sierpinski} we use \textsc{Matlab}'s built-in function \texttt{randi}, which picks from the set $\{0, 1, 2\}$. 

To calculate the new point (shown in Listing \ref{alg:triangle}), we take the average between the previous point and the corresponding vertex. For example, if the first random number generated was $2$, then we would take the average between our initial point (0, 0) and the vertex that corresponds to $2$, which is in our case is $(1, 1)$. Our new point would be $(0.5, 0.5)$. We then would repeat this process until one decides to stop.

\begin{listing}[ht]
\caption{MATLAB Code for CGR with three vertices using a random number generator}
\label{alg:triangle}
\begin{Verbatim}[commandchars=\\\{\}]
\PYG{n}{n} \PYG{p}{=} \PYG{l+m+mi}{50000}\PYG{p}{;}
\PYG{n}{pts} \PYG{p}{=} \PYG{n+nb}{zeros}\PYG{p}{(}\PYG{l+m+mi}{2}\PYG{p}{,} \PYG{n}{n}\PYG{p}{);}
\PYG{n}{rseq} \PYG{p}{=} \PYG{n}{randi}\PYG{p}{([}\PYG{l+m+mi}{0}\PYG{p}{,} \PYG{l+m+mi}{2}\PYG{p}{],} \PYG{l+m+mi}{1}\PYG{p}{,} \PYG{n}{n}\PYG{o}{\PYGZhy{}}\PYG{l+m+mi}{1}\PYG{p}{);}
\PYG{k}{for} \PYG{n+nb}{i} \PYG{p}{=} \PYG{l+m+mi}{2}\PYG{p}{:}\PYG{n}{n}
	\PYG{k}{if} \PYG{n}{rseq}\PYG{p}{(}\PYG{n+nb}{i}\PYG{o}{\PYGZhy{}}\PYG{l+m+mi}{1}\PYG{p}{)} \PYG{o}{==} \PYG{l+m+mi}{0}
		\PYG{n}{pts}\PYG{p}{(:,} \PYG{n+nb}{i}\PYG{p}{)} \PYG{p}{=} \PYG{l+m+mf}{0.5}\PYG{o}{*}\PYG{p}{(}\PYG{n}{pts}\PYG{p}{(:,} \PYG{n+nb}{i}\PYG{o}{\PYGZhy{}}\PYG{l+m+mi}{1}\PYG{p}{)} \PYG{o}{+} \PYG{p}{[}\PYG{l+m+mi}{0}\PYG{p}{;} \PYG{l+m+mi}{0}\PYG{p}{]);}
	\PYG{k}{elseif} \PYG{n}{rseq}\PYG{p}{(}\PYG{n+nb}{i}\PYG{o}{\PYGZhy{}}\PYG{l+m+mi}{1}\PYG{p}{)} \PYG{o}{==} \PYG{l+m+mi}{1}
		\PYG{n}{pts}\PYG{p}{(:,} \PYG{n+nb}{i}\PYG{p}{)} \PYG{p}{=} \PYG{l+m+mf}{0.5}\PYG{o}{*}\PYG{p}{(}\PYG{n}{pts}\PYG{p}{(:,} \PYG{n+nb}{i}\PYG{o}{\PYGZhy{}}\PYG{l+m+mi}{1}\PYG{p}{)} \PYG{o}{+} \PYG{p}{[}\PYG{l+m+mi}{0}\PYG{p}{;} \PYG{l+m+mi}{1}\PYG{p}{]);}
	\PYG{k}{else}
		\PYG{n}{pts}\PYG{p}{(:,} \PYG{n+nb}{i}\PYG{p}{)} \PYG{p}{=} \PYG{l+m+mf}{0.5}\PYG{o}{*}\PYG{p}{(}\PYG{n}{pts}\PYG{p}{(:,} \PYG{n+nb}{i}\PYG{o}{\PYGZhy{}}\PYG{l+m+mi}{1}\PYG{p}{)} \PYG{o}{+} \PYG{p}{[}\PYG{l+m+mi}{1}\PYG{p}{;} \PYG{l+m+mi}{1}\PYG{p}{]);}
	\PYG{k}{end}
\PYG{k}{end}
\PYG{n}{plot}\PYG{p}{(}\PYG{n}{pts}\PYG{p}{(}\PYG{l+m+mi}{1}\PYG{p}{,} \PYG{p}{:),} \PYG{n}{pts}\PYG{p}{(}\PYG{l+m+mi}{2}\PYG{p}{,} \PYG{p}{:),} \PYG{l+s}{\PYGZsq{}k.\PYGZsq{}}\PYG{p}{,} \PYG{l+s}{\PYGZsq{}MarkerSize\PYGZsq{}}\PYG{p}{,} \PYG{l+m+mi}{1}\PYG{p}{);}
\PYG{n}{axis}\PYG{p}{(}\PYG{l+s}{\PYGZsq{}square\PYGZsq{}}\PYG{p}{);}
\PYG{n}{set}\PYG{p}{(}\PYG{n}{gca}\PYG{p}{,} \PYG{l+s}{\PYGZsq{}xtick\PYGZsq{}}\PYG{p}{,} \PYG{p}{[]);}
\PYG{n}{set}\PYG{p}{(}\PYG{n}{gca}\PYG{p}{,} \PYG{l+s}{\PYGZsq{}ytick\PYGZsq{}}\PYG{p}{,} \PYG{p}{[]);}
\end{Verbatim}
\end{listing}

\begin{figure}[t]
	\centering
	\begin{subfigure}[t]{0.475\textwidth}
		\includegraphics[width= \textwidth]{figures/corner_serpinski_crop.png}
		\caption{Sierpinski triangle ($n = 50000$)}
		\label{fig:corner_sierpinski}
	\end{subfigure}
	\quad
	\begin{subfigure}[t]{0.475\textwidth}
		\includegraphics[width=\textwidth]{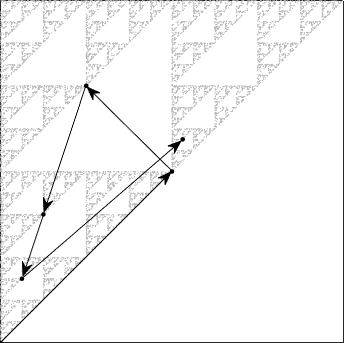}
		\caption{First 5 steps of chaos game with three vertices. In this example, the first 5 random integers generated are 2, 1, 0, 0, 2.}
		\label{fig:sierpinski_steps}
	\end{subfigure}
	\caption{CGR with three vertices of randomly generated integers.}
\end{figure}

One would expect that the outcome of Listing \ref{alg:triangle} would show that the dots would appear randomly; however, this is not the case: Figure \ref{fig:corner_sierpinski} shows the result of the chaos game if the game was played 50000 times. We can see that this figure shows a fractal-like pattern (such as one of a Sierpinski triangle), and looks almost (if not) identical to the figure which resulted from the IFS.

It is quite interesting that a sequence of random numbers could produce such a distinct pattern. Indeed, some bad \emph{pseudorandom number generators} were identified as being bad precisely because they failed a similar-in-spirit visualization test \cite{knuth1997seminumerical}. We will discuss this issue further in section~\ref{sec:random}. Figure \ref{fig:sierpinski_steps} shows an example of the first five steps of the chaos game (overlaid on top of a CGR with 10000 points), which corresponds to the following values that were randomly generated: 2, 1, 0, 0, 2. Displaying the CGR in this manner can give us some insight as to why this particular chaos game with three vertices gives us such a striking fractal: we can see from the figure that each of the points from the first five steps is located at a vertex (in this case the right-angled one) of one of the many triangles within the fractal. A careful observer would notice that the triangle at which the point is located becomes smaller each iteration. This implies that within the next few iterations, the triangle associated with the computed point would be too small to be seen in this illustration. [Note: if we zoom into that microscopic level, the dots would look random. Indeed, the random dots are ultimately dense on the Sierpinski triangle.] This sequence of points generated by the chaos game is called the \emph{orbit} of the seed and it is attracted to the Sierpinski triangle, sometimes called a \emph{strange attractor}. For another description of why this chaos game creates the Sierpinski triangle, see \cite{fisher2012science}.

\subsection{Generalization of CGR}
We have already seen that a fractal pattern can appear when we follow the chaos game described above for a polygon with three vertices (triangle) using a random number generator, but we are not limited to only this. There are many variations for CGR, though some are more useful than others. In this section, we will look into what patterns arises when we perform CGR on different polygons, changing the placement of the $m$-th point relative to the $(m-1)$-th and corresponding vertex, and explore which of these could potentially be useful for some applications presented in Section \ref{sec:DNA}.

\subsubsection{Different Polygons}
\begin{figure}
	\centering
	\begin{subfigure}[b]{0.31\textwidth}
		\includegraphics[width=\textwidth]{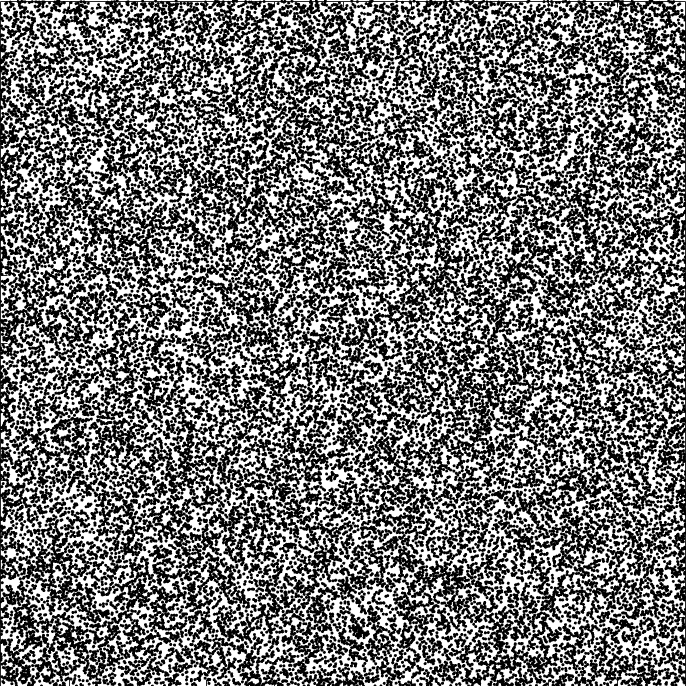}
		\caption{Square}
		\label{fig:cgr_square}
	\end{subfigure}
	~
	\begin{subfigure}[b]{0.31\textwidth}
		\includegraphics[width=\textwidth]{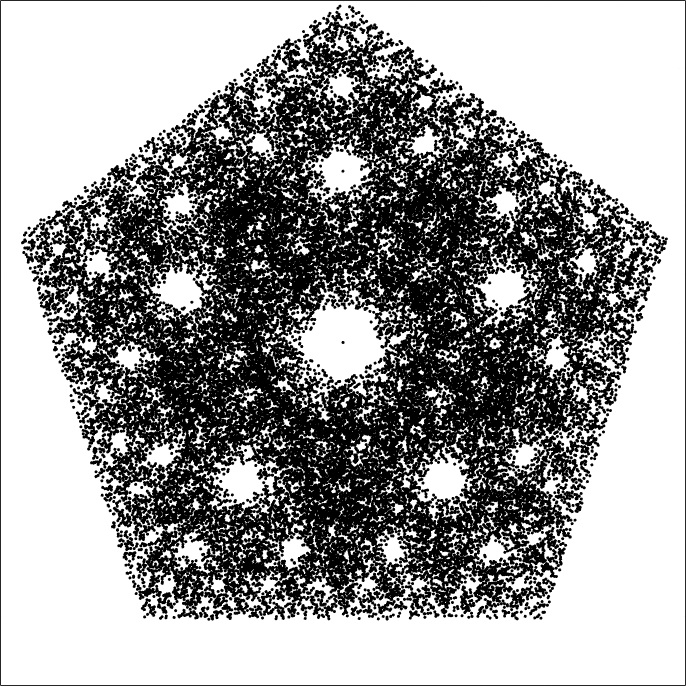}
		\caption{Pentagon}
	\end{subfigure}
	~
	\begin{subfigure}[b]{0.31\textwidth}
		\includegraphics[width=\textwidth]{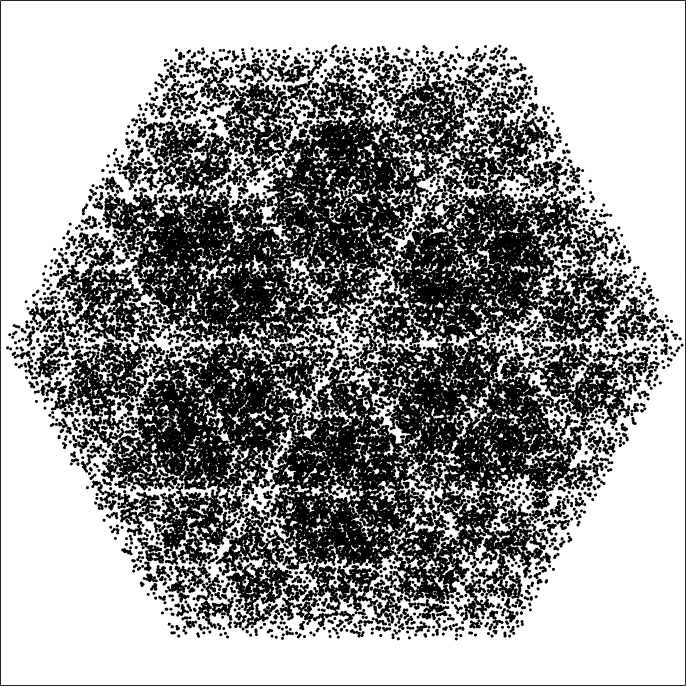}
		\caption{Hexagon}
	\end{subfigure}\\
	\begin{subfigure}[b]{0.31\textwidth}
		\includegraphics[width=\textwidth]{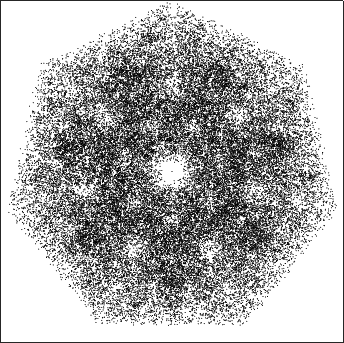}
		\caption{Heptagon}
	\end{subfigure}
	~
	\begin{subfigure}[b]{0.31\textwidth}
		\includegraphics[width=\textwidth]{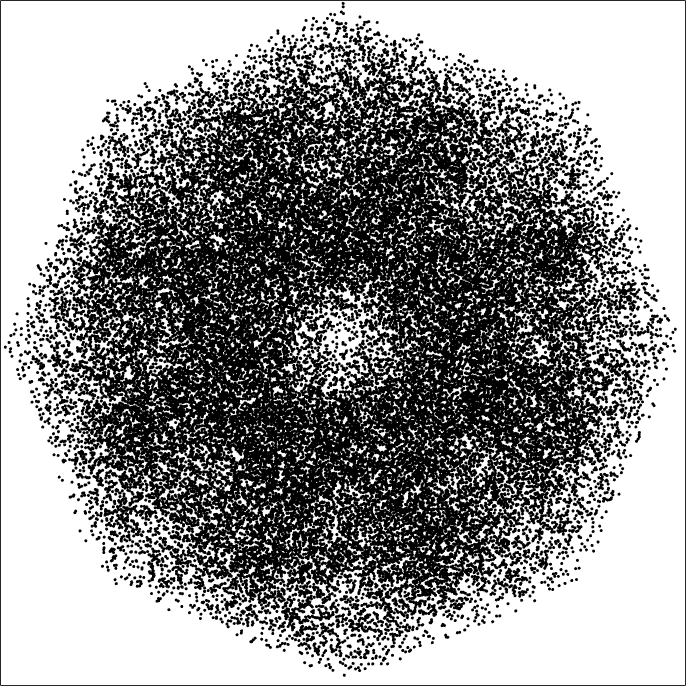}
		\caption{Octagon}
	\end{subfigure}
	~
	\begin{subfigure}[b]{0.31\textwidth}
		\includegraphics[width=\textwidth]{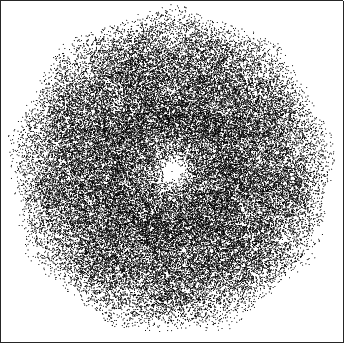}
		\caption{Nonagon}
	\end{subfigure}
	\caption{Chaos game of polygon with different number of vertices}
	\label{figs:cgr_vertices}
\end{figure}

 Figure \ref{figs:cgr_vertices} shows the chaos game representation, where the $m$-th point is placed halfway between the $(m-1)$-th point and the vertex corresponding to the $m$-th base, for various polygons, using 50,000 points. As we can see from Figure \ref{fig:cgr_square}, there is no apparent pattern: the chaos game produced a square uniformly and randomly filled with points. The other figures from Figure \ref{figs:cgr_vertices} appear to exhibit some patterns. These patterns are due to several parts of their attractors overlapping, indicated by the the uneven distribution of points (where points do appear); i.e.~we can see some darker areas and some lighter areas in the figures. Therefore, for this particular chaos game (in which we say that the dividing rate $r$ is 0.5---we will talk about this in section \ref{sec:dividingrate}), the square is the best choice to visualize one-dimensional sequences. We will see in section \ref{sec:generalization} that we can use other polygons (we will be referring to them as $n$-gons) for CGR to uncover patterns in sequences not with 4 elements by using the appropriate dividing rate.
 
\subsubsection{Dividing rate}\label{sec:dividingrate}
In the previous subsection, we saw the different patterns appear for different shapes. We saw that for shapes with more than 4 vertices had several parts of the attractor overlapping. Here we will show that attractor overlapping also happens when choosing a small dividing rate.

The chaos game is not restricted to the ``rule'' in which the new point has to be placed halfway between the current point and the corresponding vertex as depicted in the previous examples. The new point can be placed anywhere within the line segment created by the two points of reference. The placement of the new point affects how the attractor of the CGR looks. We will see this in the following example. To quantify the placement of the point, we take the proportion of which the distance between the new point and current point is from the distance between the current point and the corresponding vertex. This is called the \emph{dividing rate}, $r$. Therefore, when the new point is placed halfway between the current point and the vertex, the dividing rate is $r = 0.5$, and if the new point is placed closer to the location of the current point, then the dividing rate $r < 0.5$, whereas if the point is closer to the vertex, then the dividing rate $r > 0.5$. 

\begin{figure}
	\centering
	\begin{subfigure}[t]{0.31\textwidth}
		\includegraphics[width=\textwidth]{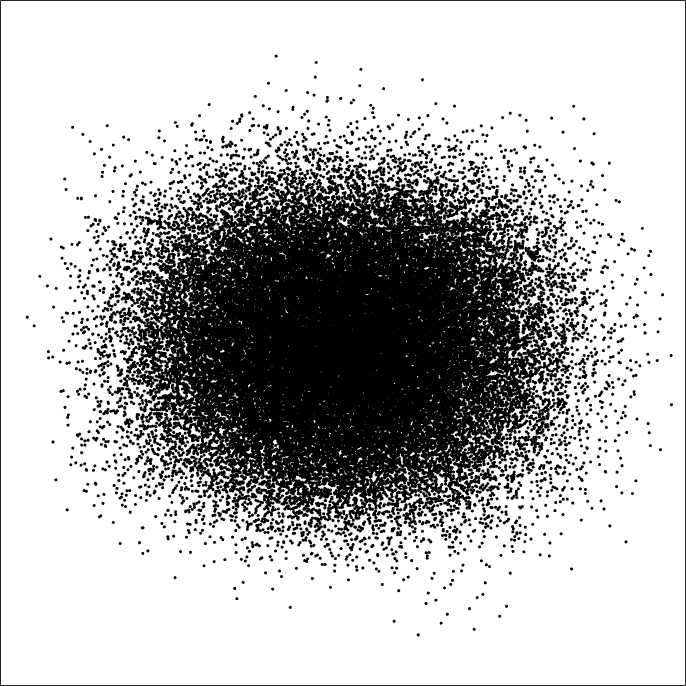}
		\caption{$r = 0.1$}
		\label{fig:square_r1}
	\end{subfigure}
	~
	\begin{subfigure}[t]{0.31\textwidth}
		\includegraphics[width=\textwidth]{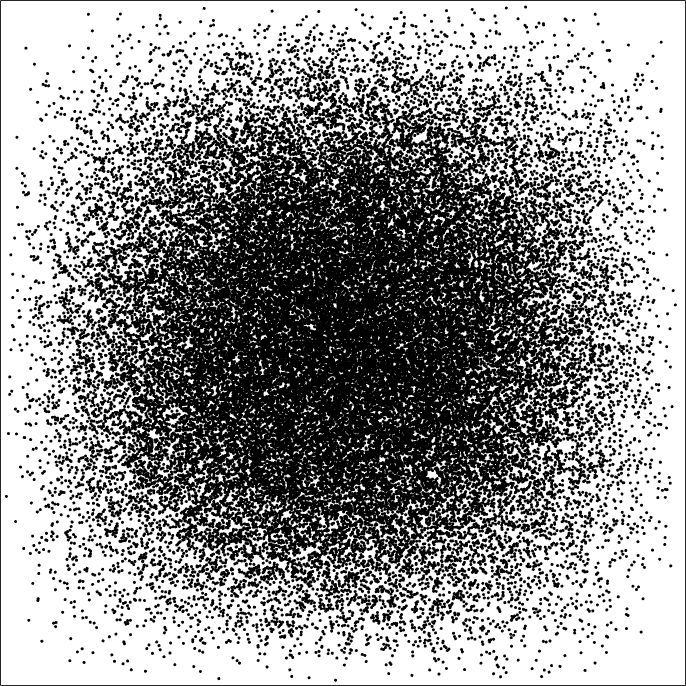}
		\caption{$r = 0.25$}
		\label{fig:square_r25}
	\end{subfigure}
	~
	\begin{subfigure}[t]{0.31\textwidth}
		\includegraphics[width=\textwidth]{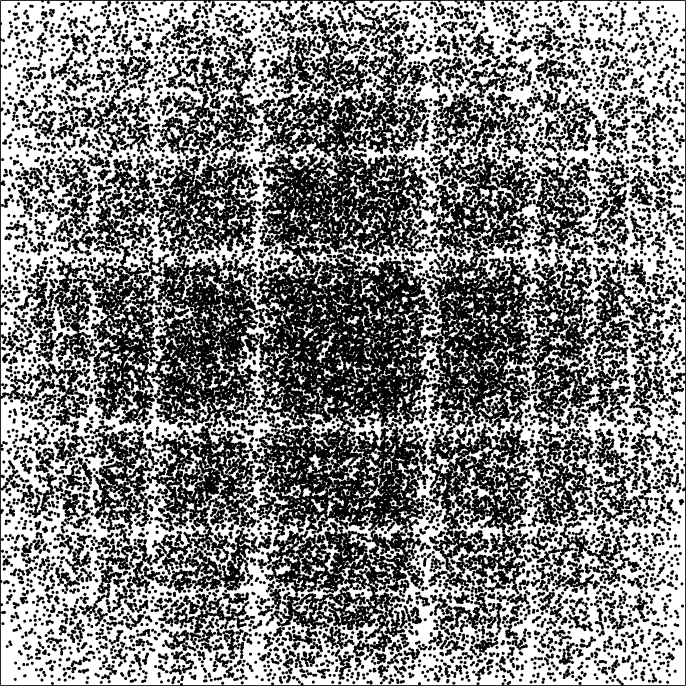}
		\caption{$r = 0.4$}
		\label{fig:square_r4}
	\end{subfigure}
	
	\begin{subfigure}[t]{0.31\textwidth}
		\includegraphics[width=\textwidth]{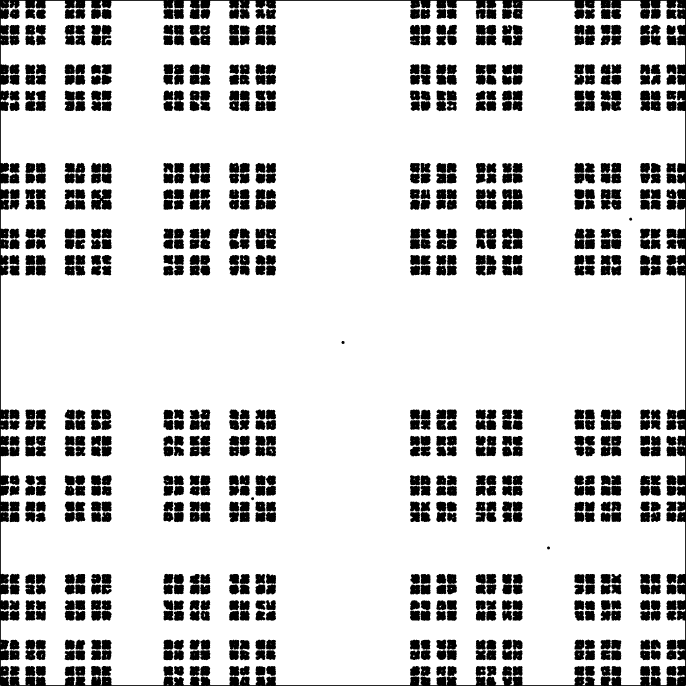}
		\caption{$r = 0.6$}
		\label{fig:square_r6}
	\end{subfigure}
	~
	\begin{subfigure}[t]{0.31\textwidth}
		\includegraphics[width=\textwidth]{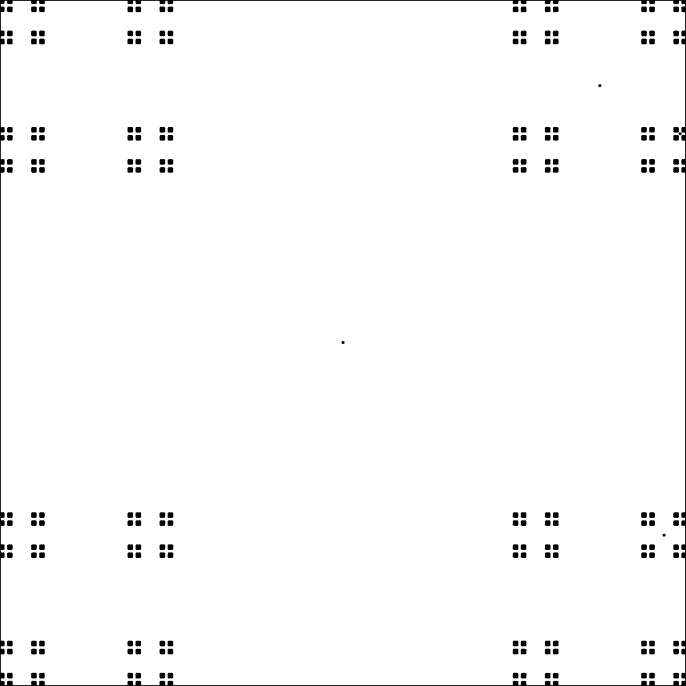}
		\caption{$r = 0.75$}
		\label{fig:square_r75}
	\end{subfigure}
	~
	\begin{subfigure}[t]{0.31\textwidth}
		\includegraphics[width=\textwidth]{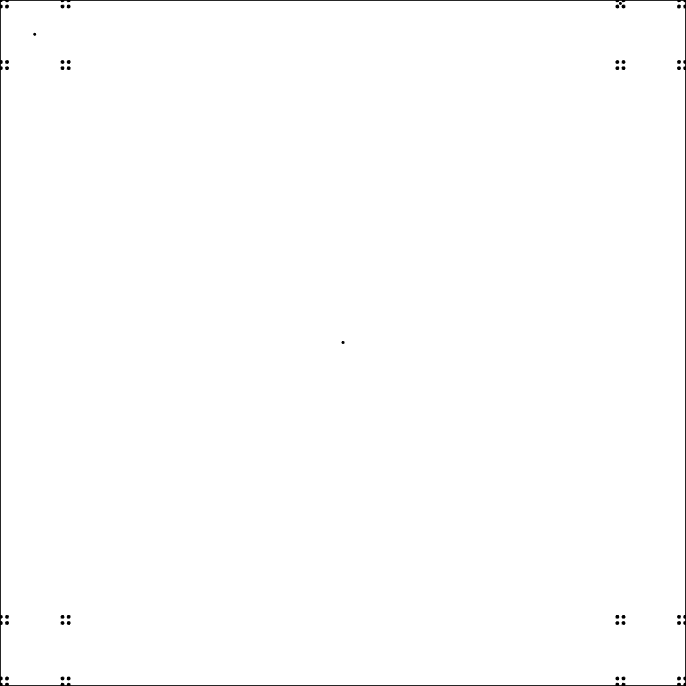}
		\caption{$r = 0.9$}
		\label{fig:square_r9}
	\end{subfigure}
	\caption{Illustration of different $r$ affecting how a square CGR looks using 50000 points.}
	\label{fig:CGR_square_R}
\end{figure}

To begin this discussion, let us vary $r$ for the square-shaped CGR of a sequence of random integers. As we have seen in the previous subsection, the square-shaped CGR did not exhibit any fractal patterns when $r = 0.5$, so we wondered if we varied $r$, would this change? Figure \ref{fig:CGR_square_R} shows the square CGRs for $r = 0.1,\ 0.25,\ 0.4,\ 0.6,\ 0.75,$ and $0.9$. From these figures, we see how~$r$ affects the visualization. 

For $r < 0.5$, we can see that some pattern appear, especially those with a dividing rate close to $r = 0.5$. We can see for $r = 0.6$, the points cover the entire plot area, but they are unevenly-space (in comparison to what we have seen for $r = 0.5$). This suggests that for this dividing rate, the attractors are overlapping, indicated by the darker regions of the plot due to the points being more densely packed. As $r$ decreases, we can see that the CGR becomes more circular with the points spaced more closely together. From this, it infers that as $r$ decreases, more of the attractors overlap.

For $r > 0.5$, we can see that a different pattern appears (compared to those with a dividing ratio $r < 0.5$). The pattern shown for the CGRs with a dividing ratio $r > 0.5$ is similar to one another, as seen in Figure \ref{fig:CGR_square_R}. Figure \ref{fig:square_r9} looks like a smaller version of Figure \ref{fig:square_r75}, which similarly, looks like a smaller version of Figure \ref{fig:square_r6}. In contrast to the CGR with the dividing rate $r > 0.5$, we can see that the points are spread apart from each other. This can possibly indicate that the attractors do not overlap, as well as do not touch each other. Although this can be used to distinguish any nonrandomness, it would be much more difficult and using the dividing rate $r = 0.5$.

In applications, it is optimal to find a dividing rate for different polygons that can offer maximum packing of nested, non-overlapping attractors in order to find a pattern in a given sequence. In the next subsection, we look at this generalization of CGR.

\subsubsection{Generalization of CGR}\label{sec:generalization}
\begin{figure}
	\centering
	\begin{subfigure}[b]{0.45\textwidth}
		\centering
		\begin{subfigure}[b]{0.45\textwidth}
			\includegraphics[width=\textwidth]{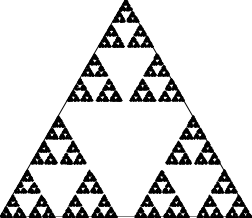} 
			\captionsetup{labelformat=empty}
			\caption{$r = 0.535898$}
		\end{subfigure}
		\quad
		\begin{subfigure}[b]{0.45\textwidth}
			\includegraphics[width=\textwidth]{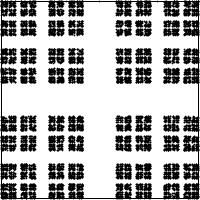} 
			\captionsetup{labelformat=empty}
			\caption{$r = 0.585786$}
		\end{subfigure}
		
		\begin{subfigure}[b]{0.45\textwidth}
			\includegraphics[width=\textwidth]{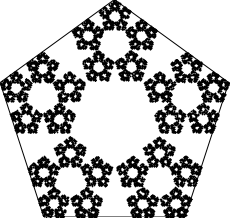} 
			\captionsetup{labelformat=empty}
			\caption{$r = 0.629808$}
		\end{subfigure}
		\quad
		\begin{subfigure}[b]{0.45\textwidth}
			\includegraphics[width=\textwidth]{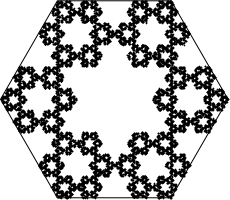} 
			\captionsetup{labelformat=empty}
			\caption{$r = 0.666666$}
		\end{subfigure}
		
		\begin{subfigure}[b]{0.45\textwidth}
			\includegraphics[width=\textwidth]{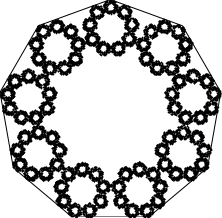} 
			\captionsetup{labelformat=empty}
			\caption{$r = 0.745145$}
		\end{subfigure}
		\quad
		\begin{subfigure}[b]{0.45\textwidth}
			\includegraphics[width=\textwidth]{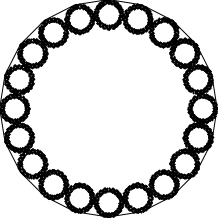} 
			\captionsetup{labelformat=empty}
			\caption{$r = 0.864726$}
		\end{subfigure}
	\end{subfigure}
	\rulesep
	\begin{subfigure}[b]{0.45\textwidth}
		\centering
		\begin{subfigure}[b]{0.45\textwidth}
			\includegraphics[width=\textwidth]{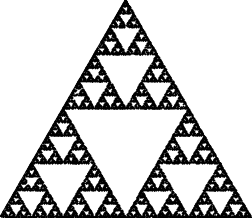} 
			\captionsetup{labelformat=empty}
			\caption{$r = 0.5$}
		\end{subfigure}
		\quad
		\begin{subfigure}[b]{0.45\textwidth}
			\includegraphics[width=\textwidth]{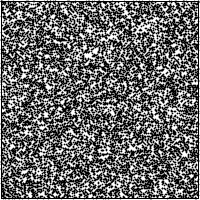} 
			\captionsetup{labelformat=empty}
			\caption{$r = 0.499999$}
		\end{subfigure}
		
		\begin{subfigure}[b]{0.45\textwidth}
			\includegraphics[width=\textwidth]{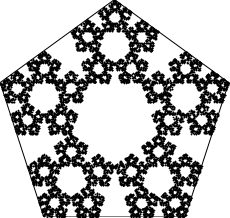} 
			\captionsetup{labelformat=empty}
			\caption{$r = 0.618033$}
		\end{subfigure}
		\quad
		\begin{subfigure}[b]{0.45\textwidth}
			\includegraphics[width=\textwidth]{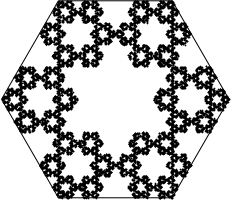} 
			\captionsetup{labelformat=empty}
			\caption{$r = 0.666666$}
		\end{subfigure}
		
		\begin{subfigure}[b]{0.45\textwidth}
			\includegraphics[width=\textwidth]{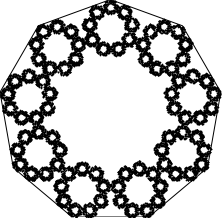} 
			\captionsetup{labelformat=empty}
			\caption{$r = 0.742227$}
		\end{subfigure}
		\quad
		\begin{subfigure}[b]{0.45\textwidth}
			\includegraphics[width=\textwidth]{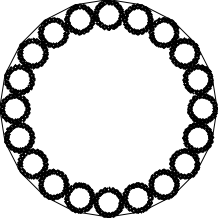} 
			\captionsetup{labelformat=empty}
			\caption{$r = 0.863271$}
		\end{subfigure}
	\end{subfigure}
	\caption{Graphical comparison of proposed CGR generalization using 50000 points. The left uses the dividing ratio formula from \cite{fiser1994chaos}, and the right uses the dividing ratio formula from \cite{almeida2009biological}.}
	\label{fig:attractor_ratio}
\end{figure}

In \cite{fiser1994chaos}, Fiser et al.~generalized CGR to be applicable for sequences of any number of elements. There are many ways of generalizing CGR, which are highlighted in \cite{almeida2009biological}, but the generalization the authors had chosen was to use an $n$-sided polygon (which they refer to as an $n$-gon), where $n$ was the number of elements in a sequence that should be represented~\cite{fiser1994chaos}. As we saw in Figure \ref{figs:cgr_vertices}, the attractors of the $n$-gons for when $n > 4$ overlapped shown by the uneven distribution of the points. Because of this, the attractor is ambiguous and is not useful to serve as a way to identity patterns in sequences. To combat this, Fiser et al.~introduced a formula to calculate the dividing rate for $n$-gons for different values of $n$:
\begin{equation*}
    r = \left(1 + \sin\left(\frac{\pi}{n}\right)\right)^{-1}\>.
\end{equation*}
Although the dividing ratio calculated from this formula does prevent attractors from overlapping, Almeida and Vinga noticed that the attractors are not optimally packed in \cite{almeida2009biological}. For example, for $n = 4$, the solution gives a dividing ratio $r = 0.585786$ instead of the \emph{typical} $r = 0.5$. In \cite{almeida2009biological}, the authors improved upon this formula for the dividing ratio:
\begin{equation*}
    r = \dfrac{2\cos\left(\pi\left(\frac{1}{2} - \frac{k}{n}\right)\right) - 2\cos\left(\pi\left(\frac{1}{2} - \frac{1}{2n}\right)\right)\cos\left((2k-1)\frac{\pi}{2n}\right)\left(1 + \frac{\tan\left((2k-1)\frac{\pi}{2n}\right)}{\tan\left(\pi-(n+2k-2)\left(\frac{\pi}{2n}\right)\right)}\right)}{2\cos\left(\pi\left(\frac{1}{2}-\frac{k}{n}\right)\right)}
\end{equation*}
where
\begin{equation*}
    k = \mathrm{round}\left(\frac{n+1}{4}\right) \>.
\end{equation*}
Figure \ref{fig:attractor_ratio} shows the generalized CGRs for various $n$-gons using the division formula from \cite{fiser1994chaos} (left) and \cite{almeida2009biological} (right).

In the next section, we will look at applications of CGR, mainly biological and mathematical. For the following examples, unless otherwise specified, assume that $r = 0.5$.


\section{CGR of biological sequences} \label{sec:DNA}
The chaos game representation has been applied to biological areas. In this section, we will look at how the chaos game has been applied the DNA and protein (amino acids) sequencing and visualizations. We also will give a brief overview of the \emph{Map of Life}, an extension of the visualization of DNA sequences. The Map of Life shows potential in allowing researchers to quantitatively classify species that was once unclear in its taxonomy.

Life certainly has random elements, but it has deterministic elements as well. The hope of CGR is to reveal some of the latter amidst the clutter of the former.

\subsection{Nucleotide sequences}
H.~J.~Jeffrey was the first to propose using the chaos game representation as a novel way of visualizing nucleotide sequences. This revealed previously unknown patterns in certain proteins~\cite{jeffrey1990chaos, jeffrey1992chaos}. The genetic sequence is made up of four bases: adenine (A), guanine (G), cytosine (C), and either thymine (T) or uracil (U) for DNA or RNA, respectively. Using a square-shaped CGR, we label the four corners with the name of each base. In this paper, A is in the bottom-left corner, C in the top-left, G in the top-right, and U/T in the bottom right. Rather than use a random number generator to determine which map to use for each iteration as we did in Section \ref{sec:chaos}, we follow the genetic sequence that we want to create a CGR for.

\begin{figure}
	\centering
	\includegraphics[width=0.8\textwidth]{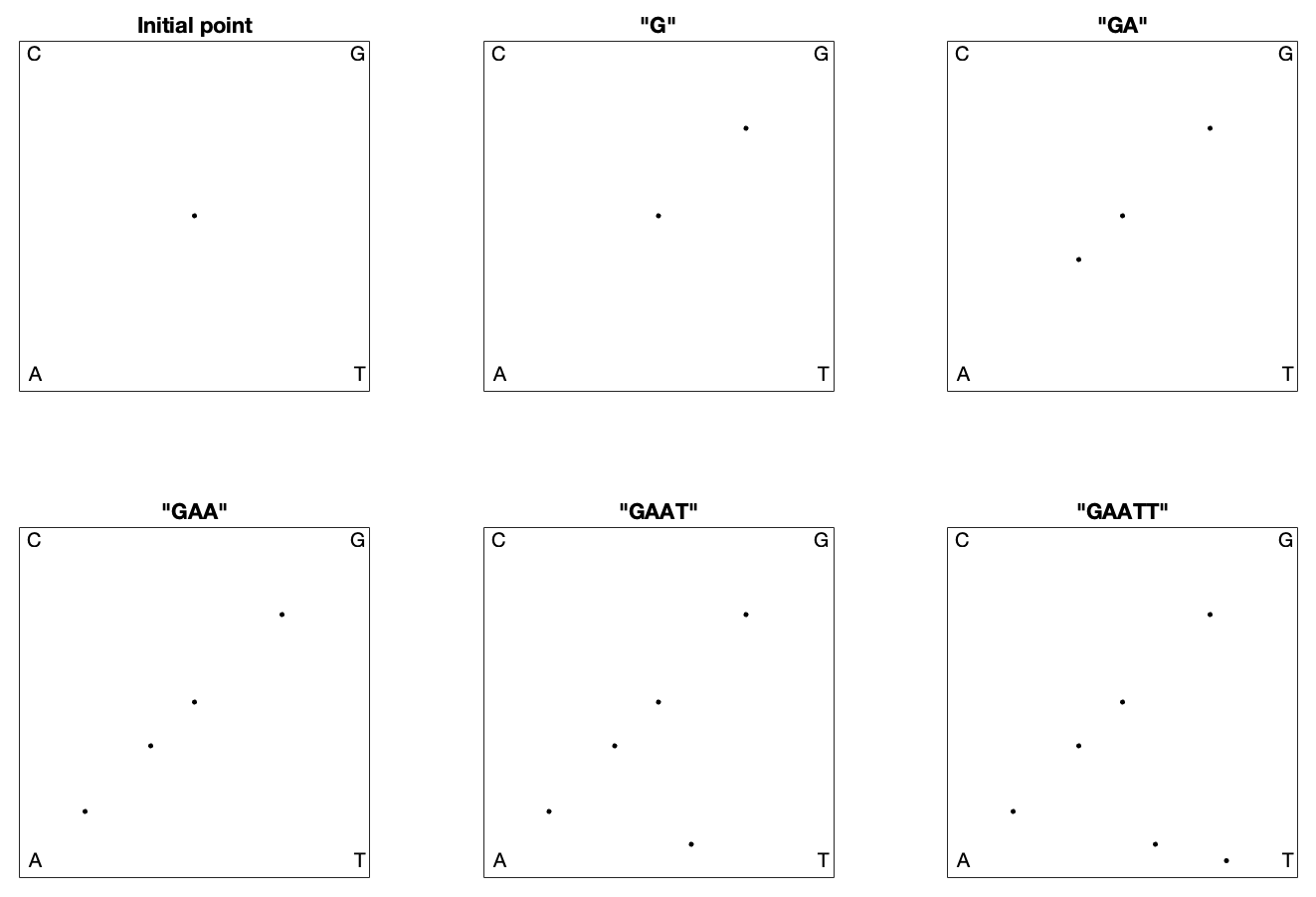}
	\caption{First five steps of producing CGR using HUMHBB DNA}
	\label{fig:intro}
\end{figure}

To demonstrate the chaos game for DNA sequences, let us walk through plotting the first five bases of the DNA sequence HUMHBB (human beta globin region, chromosome 11), ``GAATT,'' shown in Figure \ref{fig:intro}. We first mark the center as our initial point. The first base in the sequence is ``G'' so, we plot a point half way between our initial point and the ``G'' corner. The next base in the sequence is ``A'' so we plot a point half way between the point that we just plotted and the ``A'' corner. We leave the reader to check the placement of the remaining steps.

The finished CGR is shown in Figure \ref{fig:cgr_globin}, which certainly has an identifiable fractal pattern. The most prominent pattern is what is referred to as the ``double scoop,'' which actually appears in almost all vertebrate DNA sequences. This pattern is due to the fact that there is a comparative sparseness of guanine following cytosine in the gene sequence since CG dinucleotides are prone to methylation and subsequently mutation. 

\begin{figure}
	\centering
	\begin{subfigure}[t]{0.475\textwidth}
		\includegraphics[width=\textwidth]{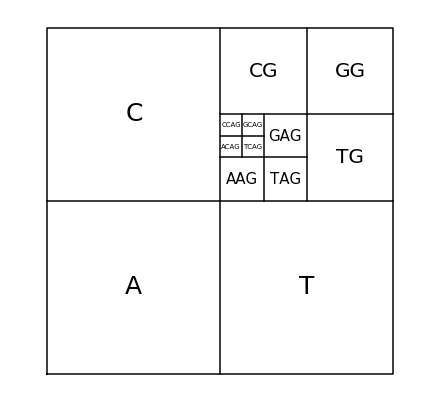}
		\caption{Correspondence between DNA\ sequences and areas of the CGR of DNA sequences}
		\label{fig:cgr_quadrants}
	\end{subfigure}
	\begin{subfigure}[t]{0.475\textwidth}
		\includegraphics[width=\textwidth]{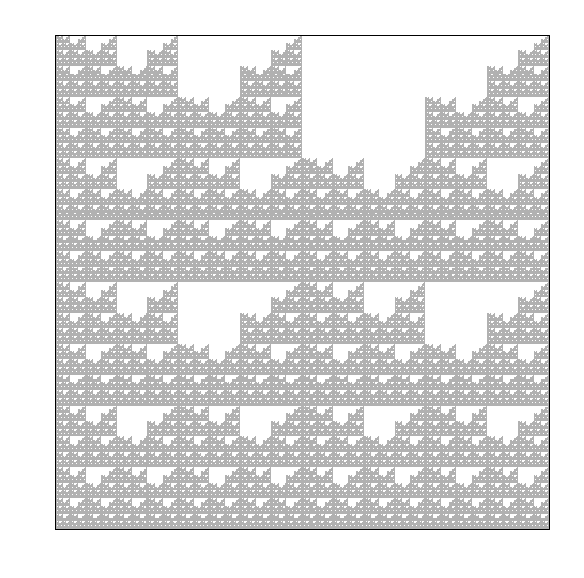}
		\caption{Explanation of the double scoop pattern---plot of CGR square with all CG quadrants unfilled.}
		\label{fig:cgr_no_cg}
	\end{subfigure}
	\caption{CGR quadrants to explain biological phenomenon such as double scoop pattern}
\end{figure}

To fully understand the double scoop pattern, we must understand the biological meaning of the CGR. Each point plotted in the CGR corresponds to a base, and depending on where it is placed, we can trace back and figure out parts of the sequence we are examining~\cite{goldman1993nucleotide}. Figure \ref{fig:cgr_quadrants} shows the relationship corresponding area of the CGR and the DNA sequence. In reference to this figure, we can see that for any point that corresponds with base G, it will be located in the upper right quadrant of the CGR plot. To see what the previous base is, we can divide the quadrant into sub-quadrants (labeling them in the same order as the quadrants), and depending on where the point is, we can determine what the previous base of the sequence is. We can repeat this step again and again to find the order in which the bases appear in the sequence.  Figure \ref{fig:cgr_no_cg} shows a CGR square where all the CG quadrants are unfilled. We can see from this figure that even though it is not a real CGR (since we did not use a sequence to produce it), we still get the same double scoop pattern found in Figure \ref{fig:cgr_globin}. By identifying regions of the CGR square in this way, it is possible to identify features of DNA sequences that correspond to patterns of the CGR.

\begin{figure}[t]
	\centering
    \begin{subfigure}[t]{0.475\textwidth}
		\includegraphics[width=\textwidth]{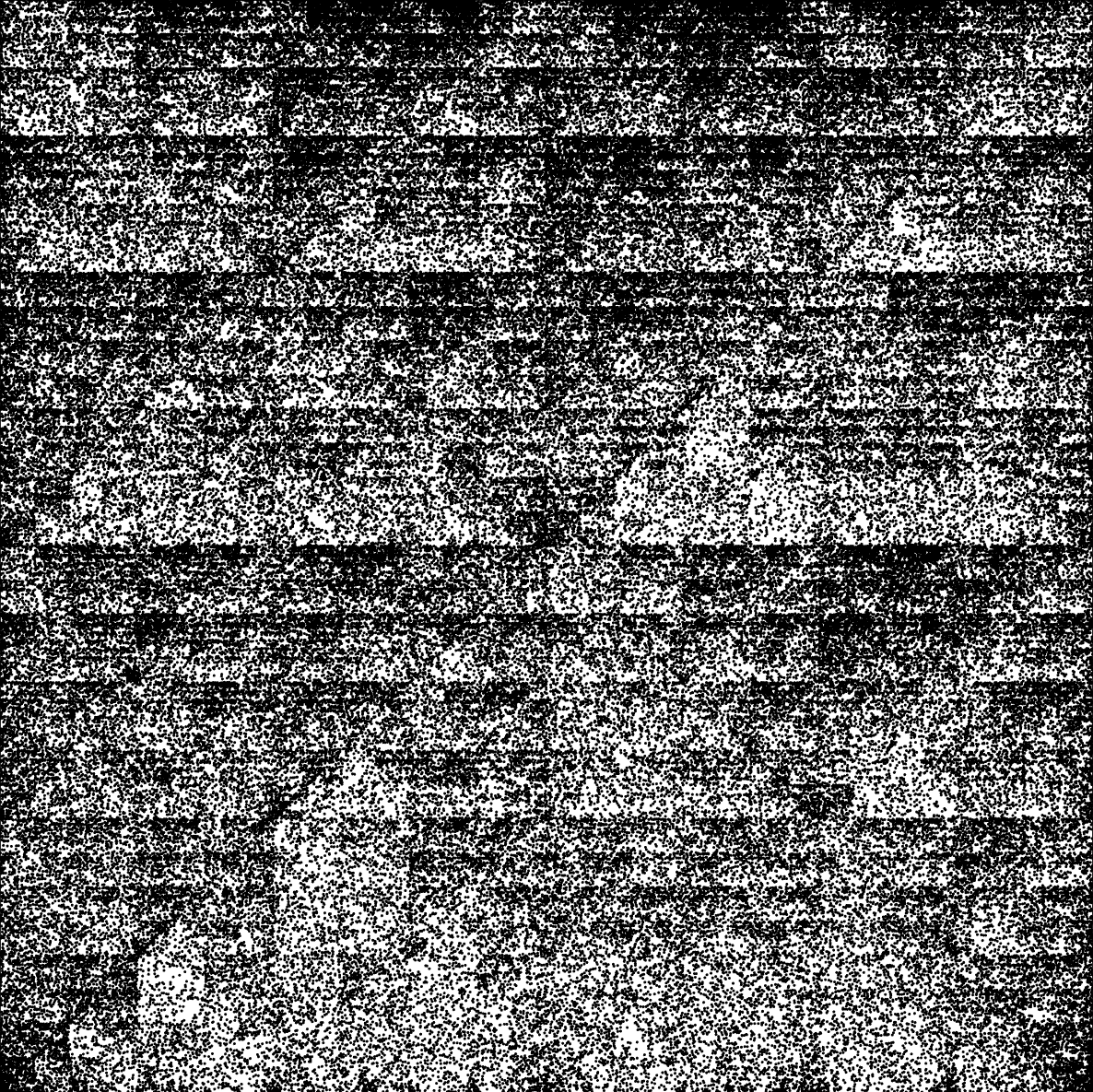}
		\caption{DNA of Human herpesvirus 5 strain TR, complete genome---235681 bp}
		\label{fig:cgr_herpesvirus}
	\end{subfigure}
    \quad
    \begin{subfigure}[t]{0.475\textwidth}
        \includegraphics[width=\textwidth]{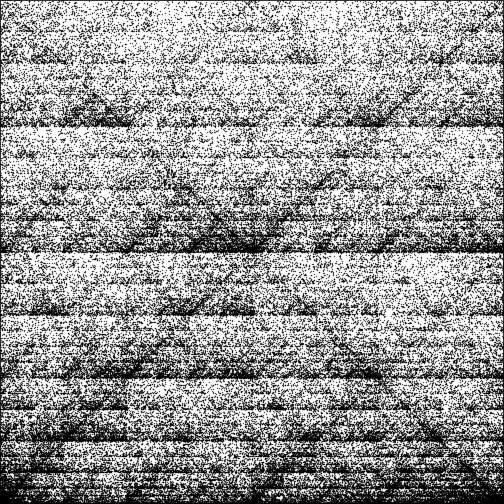}
        \caption{DNA of Chenopodium quinoa cultivar Real Blanca chloroplast, complete sequence, whole genome shotgun sequence--152282 bp}
        \label{fig:cgr_chloroplast}
    \end{subfigure}
	\caption{Examples of CGR of DNA sequences. From the bottom left corner, going clockwise, the bases are A, C, G, T.}
	\label{fig:DNA_CGR}
\end{figure}

Figure \ref{fig:DNA_CGR} shows other examples of CGRs for DNA sequences. We can see that Figure \ref{fig:cgr_herpesvirus} (CGR of DNA sequence of human herpesvirus strain) also exhibit the double scoop pattern that we have seen for HUMHBB. This agrees with what was mentioned earlier in that it is common to see the lack of C and G dinucleotides together within the human genome. Figure \ref{fig:cgr_chloroplast}, on the other hand, shows the CGR of the DNA sequence of the chloroplast of quinoa plant, and does not show a double scoop pattern. We can see from using the CGR of DNA sequences, we are able to distinguish between different species very easily.

\subsection{Map of Life}
Taking DNA sequence visualization one step further, the authors of \cite{kari2013map} proposed a novel combination of methods to create \emph{Molecular Distance Maps}. Molecular distance maps visually illustrate the quantitative relationships and patterns of proximities among the given genetic sequences and among the species they represent. To compute and visually display relationships between DNA sequences, the three techniques that were used include chaos game representation, structural dissimilarity index (DSSIM), and multi-dimensional scaling. In the paper \cite{kari2013map}, this method was applied to a variety of cases that have been historically controversial and was there demonstrated to have the potential for taxonomical clarification. In the following, we will discuss the methods used from this paper, more particularly the structural dissimilarity index and the multi-dimensional scaling as we have already looked into the chaos game representation for nucleotides.

To understand the structural dissimilarity index, we will first have to explain what the structural similarity index is. The structural similarity (SSIM) index is a method for measuring the similarity between two images based on the computation of three terms, namely luminance distortion, contrast distortion, and linear correlation~\cite{wang2004image}. It was designed to perform similarly to the human visual system, which is highly adapted to extract structural information. The overall index is a multiplicative combination of the three terms:
\begin{equation*}
	\operatorname{SSIM}(x, y) = [l(x, y)]^\alpha\cdot[c(x, y)]^\beta\cdot[s(x, y)]^\gamma \>,
\end{equation*}
where
\begin{align}
	l(x, y) &= \dfrac{2\mu_{x}\mu_{y}+C_1}{\mu_{x}^{2} + \mu_{y}^{2} + C_1} \\
	c(x, y) &= \dfrac{2\sigma_{x}\sigma_{y}+C_{2}}{\sigma_{x}^{2} + \sigma_{y}^{2} + C_2} \\
	s(x, y) &= \dfrac{\sigma_{xy} + C_{3}}{\sigma_{x}\sigma_{y} + C_{3}} \>,
\end{align}
where $\mu_{x}$, $\mu_{y}$, $\sigma_{x}$, $\sigma_{y}$ and $\sigma_{xy}$ are the local means, standard deviations, and cross-covariance for images $x$ and $y$. $C_{1}$, $C_{2}$, and $C_{3}$ are the regularization constants for the luminance, contrast, and structural terms, respectively. If $\alpha = \beta = \gamma = 1$, and $C_{3} = \sfrac{C_{2}}{2}$, the index simplifies to:
\begin{equation*}
	\operatorname{SSIM}(x, y) = \dfrac{(2\mu_{x}\mu_{y} + C_{1})(2\sigma_{xy} + C_{2})}{(\mu_{x}^{2} + \mu_{y}^{2} + C_{1})(\sigma_{x}^{2} + \sigma_{y}^{2} + C_{2})} \>.
\end{equation*}
The theoretical range of $\operatorname{SSIM}(x, y) \in [-1, 1]$, where the high value indicates high similarity. In \cite{kari2013map}, instead of calculating the overall SSIM, they computed the local SSIM value for each pixel in the image $x$. They refer to this as a distance matrix. Both global and local SSIM index can be computed in \textsc{Matlab} by using the function \texttt{ssim}. For our own experiments presented in this paper, we use the global SSIM value rather than the local.

Now that we know how to compute the SSIM index, we are able to compute the structural dissimilarity (DSSIM) index:
\begin{equation*}
	\operatorname{DSSIM}(x, y) = 1 - \operatorname{SSIM}(x, y) \>,
\end{equation*}
whose theoretical range is $[0, 2]$ with the distance being 0 between two identical images, and 2 if the two images are negatively correlated. An example of two images that are negatively correlated would be if one is completely white, while the other is completely black. Typically, the range that DSSIM falls with is $[0, 1]$---the authors from \cite{kari2015mapping} noted that almost all (over 5 million) distances that found were between $0$ and $1$, with only half a dozen exceptions of distances between $1$ and $1.0033$.

\begin{figure}
    \centering
    \begin{subfigure}[t]{0.475\textwidth}
		\includegraphics[width=\textwidth]{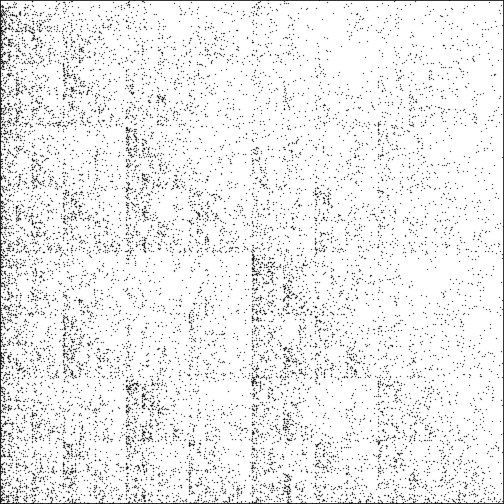}
		\caption{CGR of human (Homo sapiens sapiens) mitochondrial DNA --- 16,569 bp}
		\label{fig:cgr_human}
	\end{subfigure}
	\quad
	\begin{subfigure}[t]{0.475\textwidth}
	    \includegraphics[width=\textwidth]{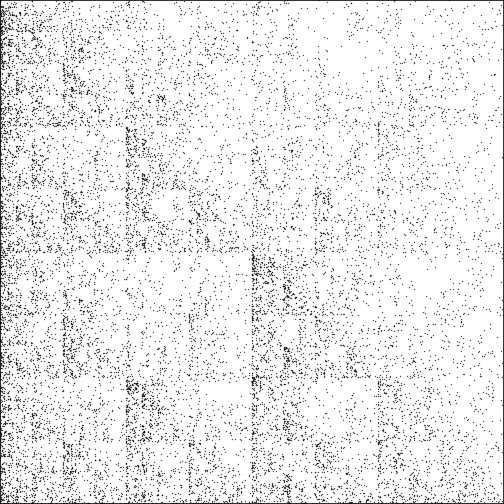}
		\caption{CGR of neanderthal (Homo sapiens neanderthalensis) --- 16,565 bp}
		\label{fig:cgr_neanderthal}
	\end{subfigure}
	\caption{CGR of human and neanderthal mitochondrial DNA.}
	\label{fig:CGR_homo_sapien}
\end{figure}

To demonstrate this, let us compare the ``anatomically modern" human (Homo sapiens sapiens) mitochondrial DNA and the neanderthal (Homo sapiens neanderthalensis) mitochondrial DNA (mtDNA). As seen from their scientific classification, we can see that these two species are from the same genus, which suggests that their DNA should be similar. Figure \ref{fig:CGR_homo_sapien} shows the CGRs of the ``anatomically modern" human (\ref{fig:cgr_human}) and of the neanderthal (\ref{fig:cgr_neanderthal}). We can see that the two CGRs are quite similar (which is agrees with the fact that the two species are from the same genus) and therefore, we expect the the DSSIM to be small. Using \textsc{Matlab}'s \texttt{ssim} function, the structural similarity index between the two images is $0.8777$, so the structural dissimilarity index is $1 - 0.8777 = 0.1223$, which agrees with our expectations.

\begin{figure}
    \centering
    \begin{subfigure}[t]{0.475\textwidth}
		\includegraphics[width=\textwidth]{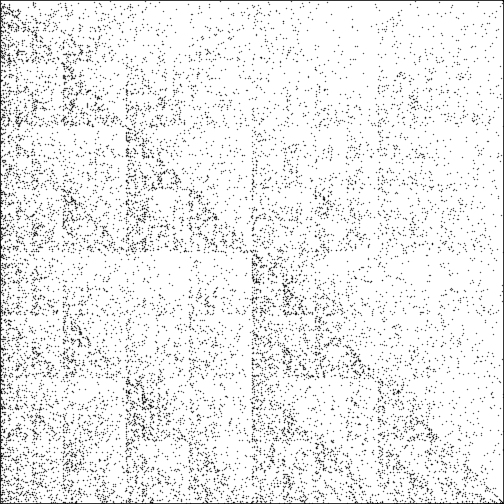}
		\caption{CGR of the great spotted kiwi (Apteryx haastii) mitochondrial DNA --- 16,980 bp}
		\label{fig:cgr_kiwi}
	\end{subfigure}
	\quad
	\begin{subfigure}[t]{0.475\textwidth}
	    \includegraphics[width=\textwidth]{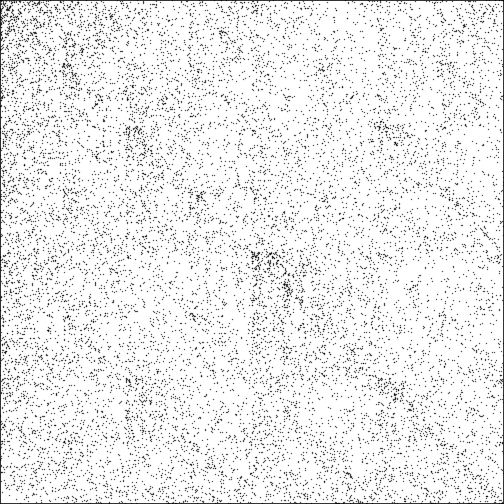}
		\caption{CGR of the pearlfish (Carapus bermudensis) --- 16,613 bp}
		\label{fig:cgr_pearlfish}
	\end{subfigure}
	\caption{CGR of other species' mitochondrial DNA.}
	\label{fig:CGR_animals}
\end{figure}

We can look at another example that compares the mtDNA of humans to other species, such as the great spotted kiwi (Figure \ref{fig:cgr_kiwi}) and pearlfish (Figure \ref{fig:cgr_pearlfish}). Just visually, the CGR of the mtDNA of the great spotted kiwi looks more similar to the CGRs of the human and neanderthal than the CGR of the pearlfish. Using all of the DSSIM, we can produce what is called a distance matrix, which is a square, real, symmetric matrix where each element is the structural dissimilarity index of two species with corresponds to the row and column. The distance matrix for the four species that we have already produced the CGR for is
\begin{equation*}
    \mathrm{DSSIM}(x, y) =
    \begin{bmatrix}
        0 & 0.1223 & 0.7823 & 0.8541 \\
        0.1223 & 0 & 0.7821 & 0.8533 \\
        0.7823 & 0.7821 & 0 & 0.8519 \\
        0.8541 & 0.8533 & 0.8519 & 0
    \end{bmatrix} \>,
\end{equation*}
where the first row/column represents humans, the second row/column represents the neanderthals, the third row/column represents the great spotted kiwi, and lastly, the fourth row/column represents the pearlfish. Obviously, the diagonal of the matrix would be $0$, since the DSSIM of two of the same images would produce a result of 0. As we have seen previous, the $\mathrm{DSSIM}(\mathrm{human},\ \mathrm{neanderthal}) = 0.1223$, and vice versa. From this matrix, we can see that the CGR of the mtDNA of pearlfish is most different than the rest---this agrees with what we have observed.

Multi-dimensional scaling (MDS) is a means of visualizing the level of similarity or dissimilarity of individual case of a dataset. It has been used in various fields such as cognitive science, information science, pychometrics, marketing, ecology, social science and other areas of study \cite{borg2005modern}. The goal of MDS is to find a spatial configuration of objects when all that is known is some measure of their general similarity or dissimilarity \cite{wickelmaier2003introduction}. In our case, MDS takes in the distance matrix and outputs a two-dimensional map, where each item is represented by a point. In \cite{kari2013map}, the authors used a classical MDS, which assumes that all the distances (from the distance matrix) are Euclidean. For the algorithm for classical MDS, refer to \cite[p.~10]{wickelmaier2003introduction}. \textsc{Matlab} has its own built-in function that does MDS, called \texttt{cmdscale}. 

\begin{figure}[th]
	\centering
	\includegraphics[width=\textwidth]{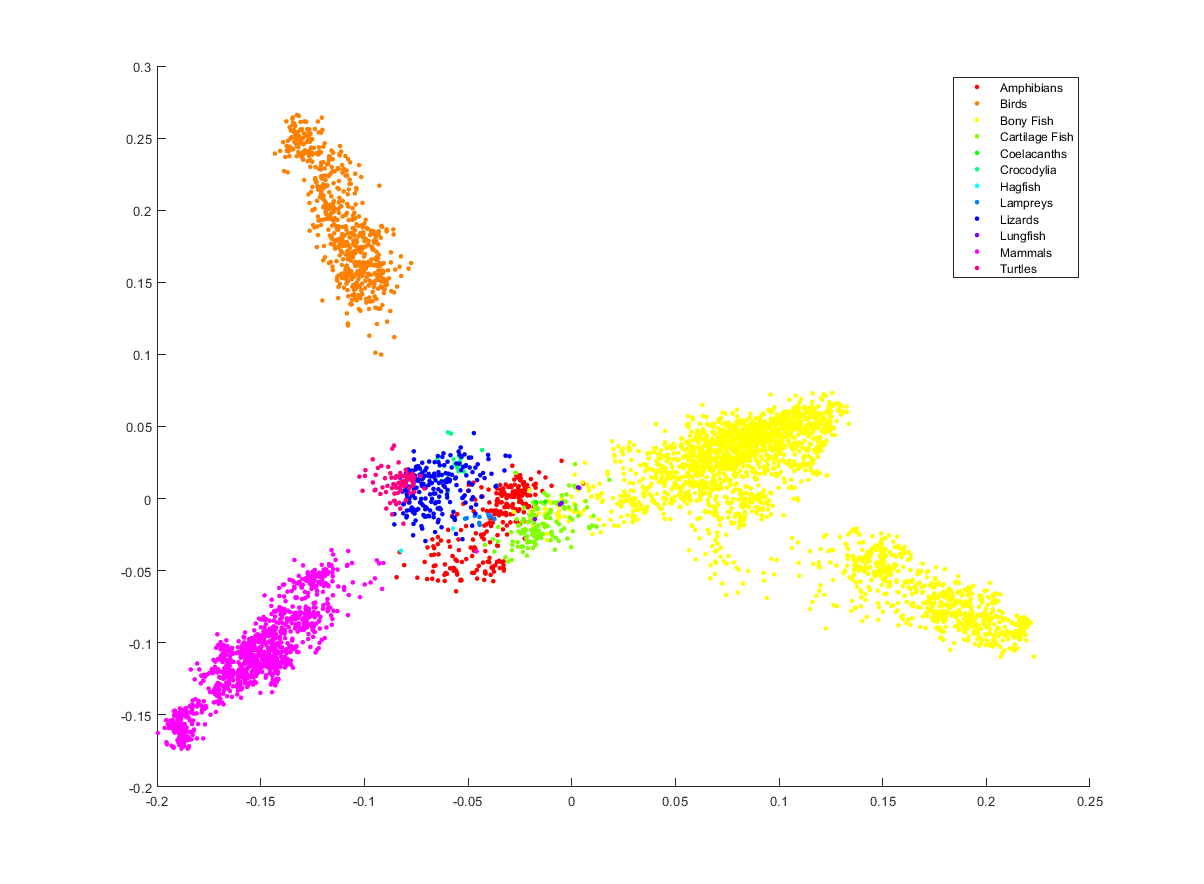}
	\caption{Molecular distance map of 4844 animals from various classes from the chordate phylum, coloured according to their class.}
	\label{fig:mdm}
\end{figure}

Figure \ref{fig:mdm} shows an example of a molecular distance map of 4844 animals from various classes from the chordate phylum (meaning that all these animals have a vertebrate), coloured according to their class. The sequences were taken from NCBI Reference Sequence Database (RefSeq). We can see from this figure that this method (Map of Life) does quite a good job in sorting species into different categories. One particular aspect that we want to point out in this figure is where the points representing the lungfish are---they are in the area of where the bony fish, cartilage fish, and amphibians touch. This is truly remarkable as lung fish has qualities of all these classes.

\subsection{Protein visualization}
\begin{figure}
	\centering
	\begin{subfigure}[t]{0.3\textwidth}
		\includegraphics[width=\textwidth]{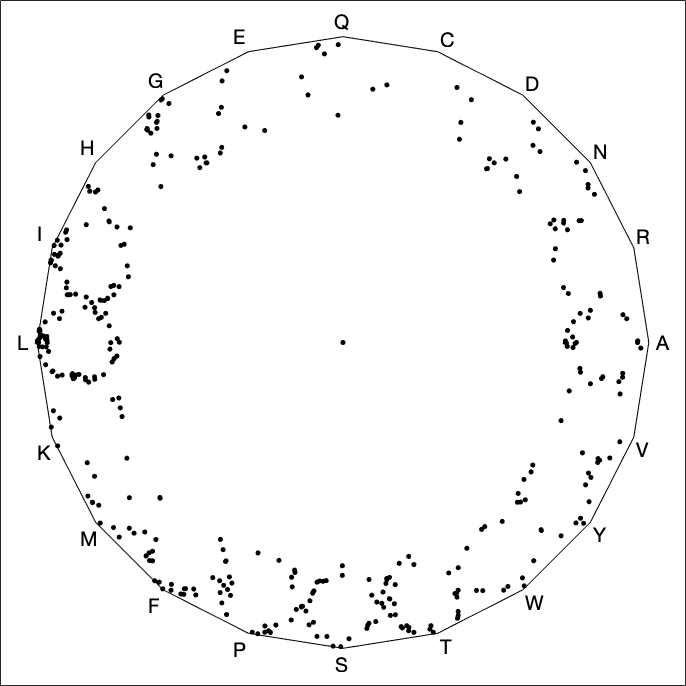}
		\caption{cytochrome b (mitochondrion) [Homo sapiens] (GenBank: ASY00349.1)}
	\end{subfigure}
	~
	\begin{subfigure}[t]{0.3\textwidth}
		\includegraphics[width=\textwidth]{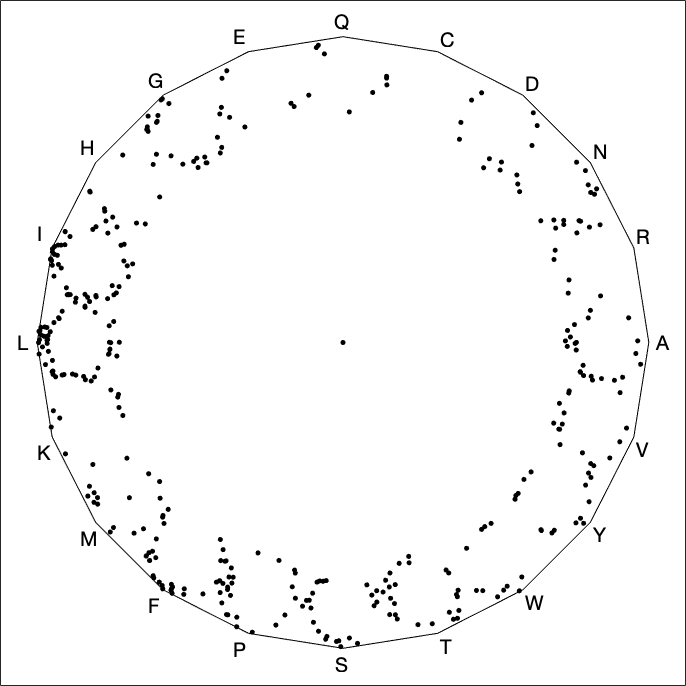}
		\caption{cytochrome b (mitochondrion) [Mus musculus] (GenBank: NP\textunderscore904340.1)}
	\end{subfigure}
	~
	\begin{subfigure}[t]{0.3\textwidth}
		\includegraphics[width=\textwidth]{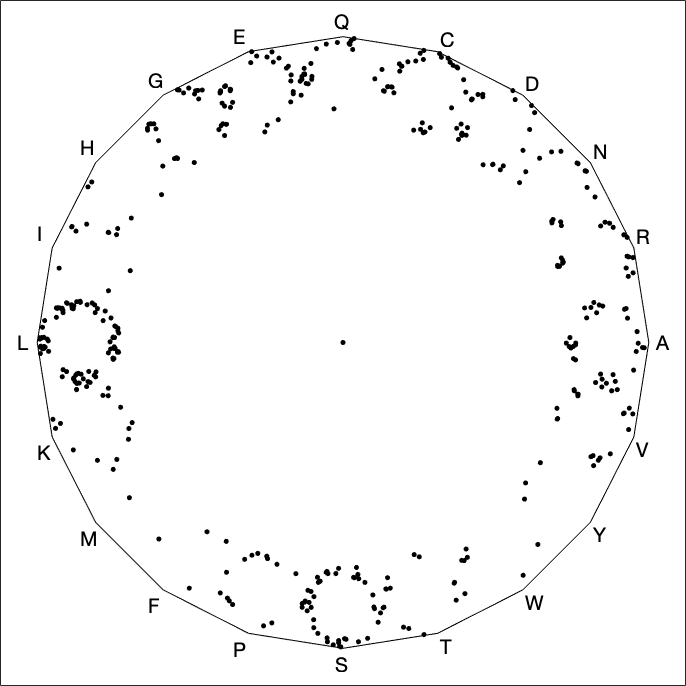}
		\caption{ribonuclease inhibitor [Homo sapien] (GenBank: NP\_976317.1)}
	\end{subfigure}
	
	\begin{subfigure}[t]{0.3\textwidth}
		\includegraphics[width=\textwidth]{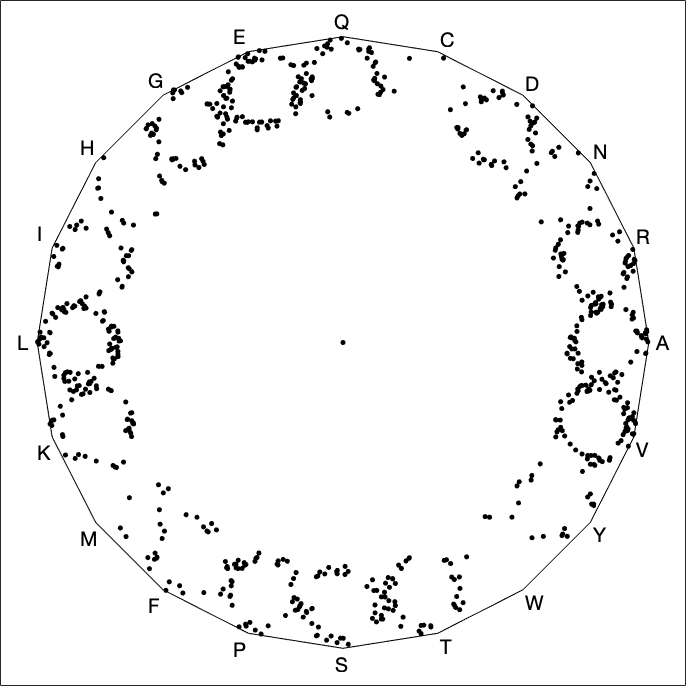}
		\caption{major vault protein [Rattus norvegicus] (GenBank: NP\_073206.2)}
	\end{subfigure}
	~
	\begin{subfigure}[t]{0.3\textwidth}
		\includegraphics[width=\textwidth]{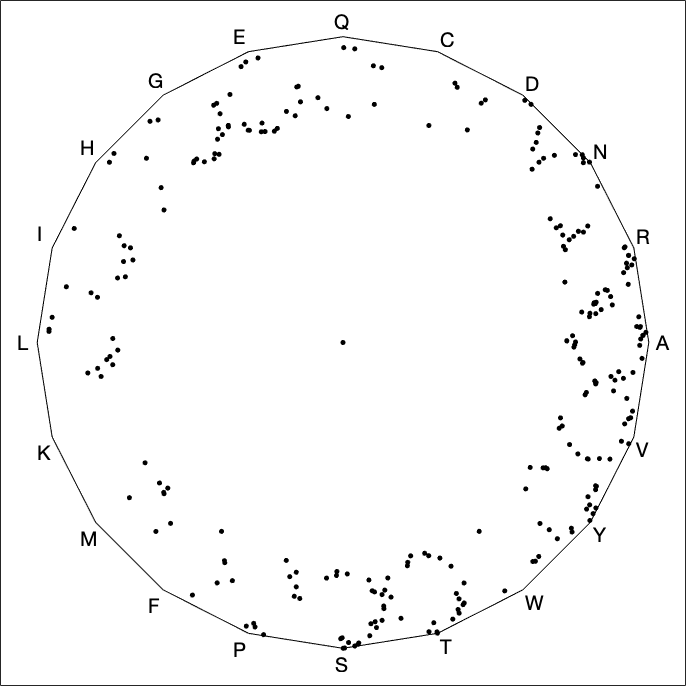}
		\caption{pertussis toxin subunit 1 [Bordetella pertussis Tohama I] (GenBank: NP\_882282.1)}
	\end{subfigure}
	~
	\begin{subfigure}[t]{0.31\textwidth}
		\includegraphics[width=\textwidth]{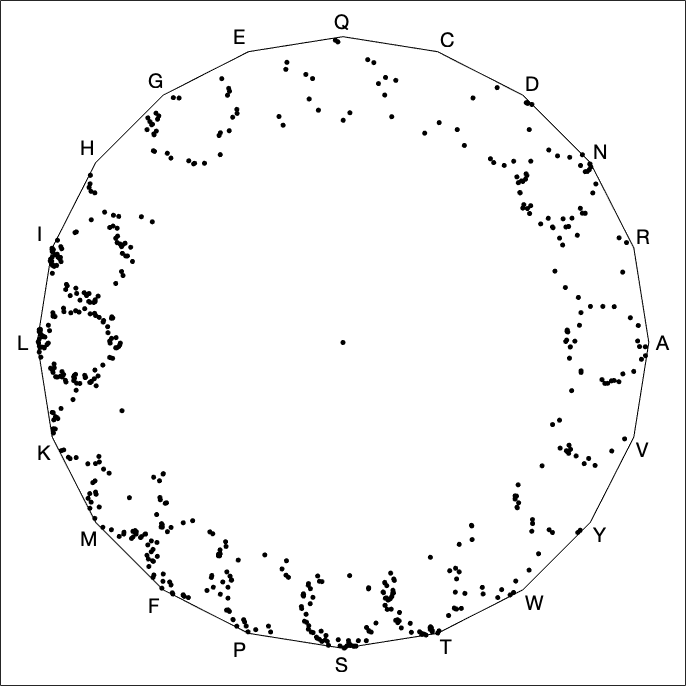}
		\caption{NADH dehydrogenase subunit 5 (mitochondrion [Mus musculus] (GenBank: NP\_904338.1)}
	\end{subfigure}
	\caption{Illustration of the chaos game representation of proteins (specifically using their amino acids).}
	\label{fig:proteinCGR}
\end{figure}

CGR has also been applied to visualizing and analyzing both the primary and secondary structures of proteins. The primary structure of a protein is simply an amino acid sequence. To analyze the primary structure of proteins, the authors from \cite{fiser1994chaos} using an $20$-sided regular polygon, each representing an amino acid. Table \ref{tab:amino_acid} shows all 20 amino acids along with their 3-letter code, as well as their 1-letter code. To avoid the attractor from overlapping itself, we will use the dividing rate for a $20$-gon shown from \cite{almeida2009biological} in Figure \ref{fig:attractor_ratio}, $r = 0.863271$. Some CGRs for protein visualizations can be found in Figure \ref{fig:proteinCGR}.

\begin{table}[]
    \centering
    \caption{Table of all 20 amino acids}
    \begin{tabular}{c|c|c}
        Name & 3 letter code & 1 letter code  \\
        \hline
        alanine & ala & A \\
        arginine & arg & R \\
        asparagine & asn & N \\
        aspartic acid & asp & D \\
        cysteine & cys & C \\
        glutamine & gln & Q \\
        glutamic acid & glu & E \\ 
        glycine & gly & G \\
        histidine & his & H \\
        isoleucine & ile & I \\
        leucine & leu & L \\
        lysine & lys & K \\
        methionine & met & M \\
        phenylalanine & phe & F \\
        proline & pro & P \\
        serine & ser & S \\
        threonine & thr & T \\
        tryptophan & trp & W \\
        tyrosine & tyr & Y \\
        valine & val & V
    \end{tabular}
    \label{tab:amino_acid}
\end{table}

According to Basu et al.~\cite{basu1997chaos}, there are two serious limitations for using a 20-vertex CGR to identify patterns. The first limitation the authors mentioned about using the 20-vertex is that it is hard to visualize different protein families since one would have to plot the sequences in different polygons rather than plotting a random sampling of proteins of different functions and origins in a single CGR~\cite{fiser1994chaos}. The second limitation is that there the amino acid residues in different positions are often replaced by similar amino acids~\cite{basu1997chaos}. This means that the 20-vertex CGR cannot be used to differentiate between similar and dissimilar residues and the visualization could be different for proteins within the same family.

To fix the second limitation mentioned in~\cite{basu1997chaos}, Fiser et al.~\cite{fiser1994chaos} introduced using CGR to study 3D structures of proteins. A chain of amino acid is what is called a polypeptide. Due to each amino acid having a specific structure, which contributes to the amino acid's properties, such as hydrophobicity or hydrogen-bonding, depending on the sequence of the amino acid chain, the polypeptide chain will fold into its lowest energy configuration~\cite{simon1991calculation}.
This is known as the secondary structure. The two most common structures that appear in this stage are the $\alpha$-helices and the $\beta$-sheets. Using the chaos game on the secondary structure can indicate any non-randomness of the structural elements in proteins. One way to achieve this is to divide the 20 amino acids into 4 groups based on different properties and assign each group to a corner. Some properties that have been considered for CGR include hydrophobicity, molecular weight, isoelectric point (pI), $\alpha$ propensity, and $\beta$ propensity. For details, see \cite{basu1997chaos, fiser1994chaos}.

\section{CGR of mathematical sequences} \label{sec:math}
In the first-year course we gave in 2015~\cite{chan2018random}, we applied the chaos game representation to mathematical sequences found on the Online Encyclopedia of Integer Sequences (OEIS)~\cite{oeis}. In the following, we show some CGR experiments that came from the first-year course. Some sequences that were used from the Online Encyclopedia of Integer Sequences (OEIS) included digits of $\pi$ (A000796), Fibonacci Numbers (A000045), prime numbers (A000040), and some from the continued fractions section (which were introduced at the beginning of the course). We will not only look at the visualizations created in the first-year class, but we will also determine whether these experiments are considered to be good visual representations. Note that for the following 4-vertex examples, the vertices are labeled from ``0'' to ``3'' starting from the bottom left corner of the plot going clockwise (i.e.~``0'' corresponds to the coordinate point $(-1, -1)$, ``1'' corresponds to the point $(-1, 1)$, etc.). 

\subsection{Digits of pi}
\begin{figure}[th]
	\centering
	\begin{subfigure}[t]{0.475\textwidth}
		\includegraphics[width=\textwidth]{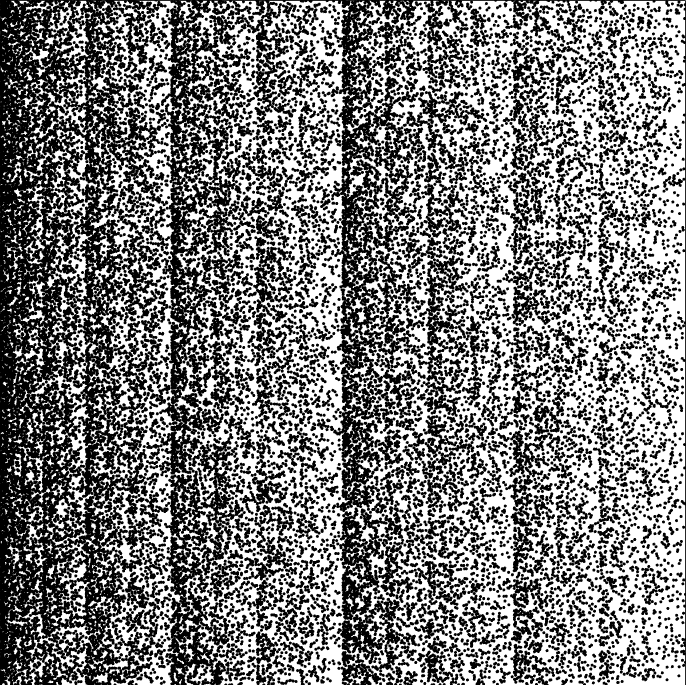}
	    \caption{4 vertex-CGR of 50000 digits of $\pi$}
	    \label{fig:digits_pi}
	\end{subfigure}
	\quad
	\begin{subfigure}[t]{0.475\textwidth}
		\includegraphics[width=\textwidth]{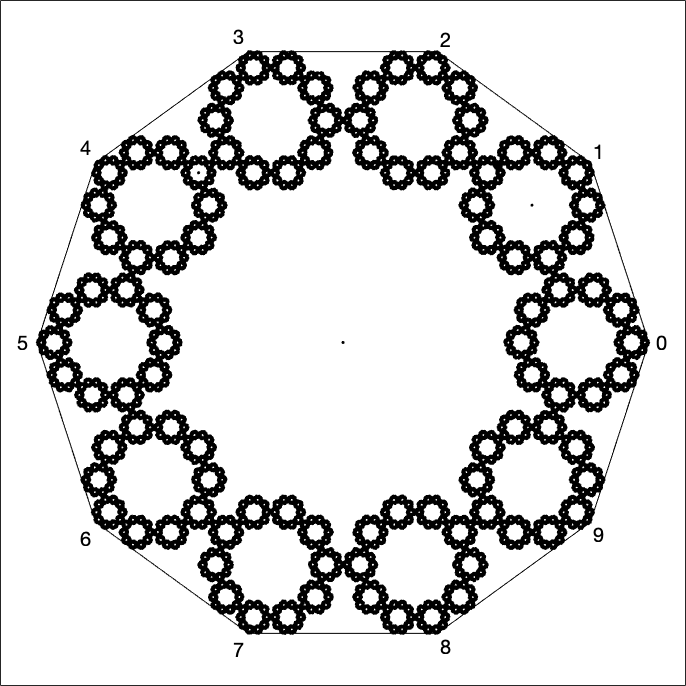}
	    \caption{10 vertex-CGR of 100000 digits of $\pi$}
	    \label{fig:digits_pi_mod10}
	\end{subfigure}
	\caption{CGR of the digits of $\pi$}
\end{figure}

We first applied CGR to the digits of $\pi$ (A000796 from the OEIS\footnote{\url{https://oeis.org/A000796}}), since we believe that the sequence of the digits are random. Since we focused on the 4-vertex CGR in the course, we na{\"i}vely took the digits of $\pi$ modulo 4 and applied it to the 4-vertex CGR. The result of this is shown in Figure \ref{fig:digits_pi}. Since we assumed that the digits of $\pi$ are random, this means that the visualization should be similar to the square-shaped CGR of random integer sequences (Figure \ref{fig:cgr_square}). However, we can see from Figure \ref{fig:digits_pi} that this is not the case; instead of a uniform covering of the space, we can see a pattern of vertical lines occurring. Quantitatively, the DSSIM index between Figures \ref{fig:digits_pi} and \ref{fig:cgr_square} is $0.9851$, which indicates that they are not similar. The reason for this is because when we take $10\bmod{4}$, the values are distributed to the corners unevenly: our 0 and 1 vertices would have one extra value each, and thus why our CGR of the digits of $\pi$ does not match the CGR of random integers even though we know for a fact that the digits of $\pi$ are indeed random. Due to this reason, Figure \ref{fig:digits_pi} is not a good visual representation of the digits of $\pi$.

Another CGR to visualize the digits of $\pi$ is shown in Figure \ref{fig:digits_pi_mod10}. We use a 10-vertex CGR ($r = 0.763932$) where each vertex represents an integer from 0 to 9. From this figure, we can see that the points produce a fractal pattern of the decagon, which suggests that the sequence is random.


\subsection{Fibonacci Sequence}
\begin{figure}
    \centering
    \includegraphics[width=0.5\textwidth]{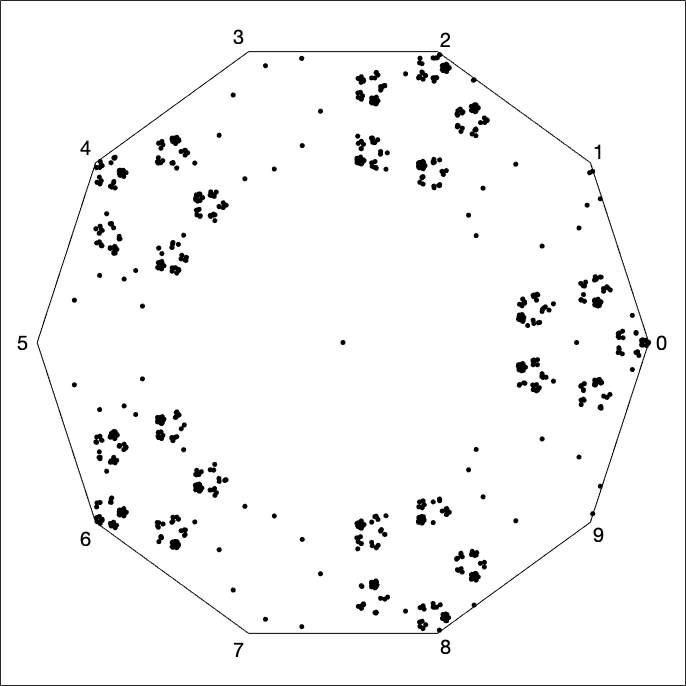}
	\caption{10 vertex-CGR of the first 3000 Fibonacci sequence modulus 10}
	\label{fig:cgr_Fibonacci}
\end{figure}

We also created a visualization of the Fibonacci Sequence (A000045 from the OEIS). One example of a visualization of this sequence is shown in Figure \ref{fig:cgr_Fibonacci}. It is a 10-vertex CGR of the Fibonacci sequence modulo 10. We can see from the figure that most of the numbers of the Fibonacci sequence are even. Other visualizations for the Fibonacci sequence can be done, such as taking the modulo of the numbers using a different divisor and using a CGR with different number of vertices.

However, since the Fibonacci sequence grows very quickly, the students learned about the limitations of floating point since they were unable to take more than the first 3000 Fibonacci numbers before the computer registering the number as infinity.

\subsection{Prime numbers}
Following the visualization of the digits of $\pi$ and the Fibonacci sequence using the chaos game, we looked at another popular sequence: the sequence of prime numbers (A000040 from the OEIS\footnote{\url{https://oeis.org/A000040}}). We looked at four different ways of creating the CGR for prime numbers. 

\begin{figure}
	\centering
	\begin{subfigure}[t]{0.475\textwidth}
		\includegraphics[width=\textwidth]{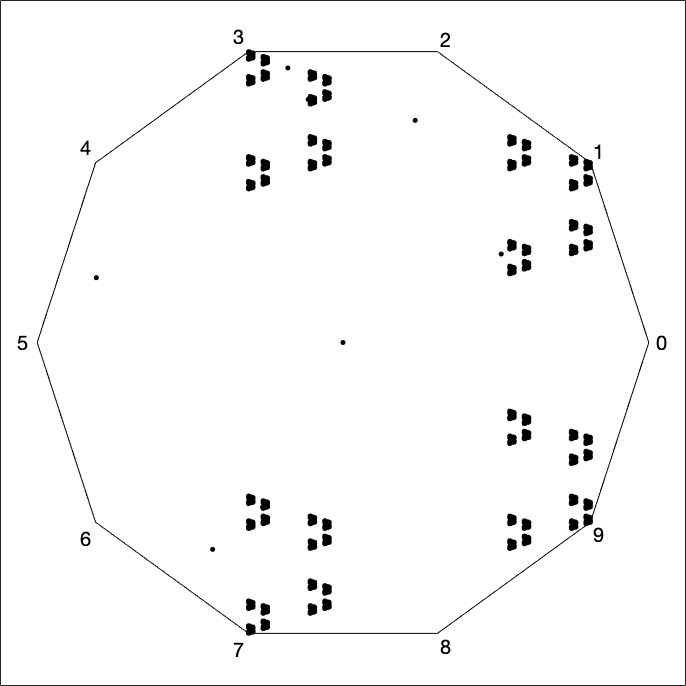}
		\caption{All prime numbers less than $1000000$ taken~$\bmod\ 10$}
		\label{fig:cgr_prime_mod_10}
	\end{subfigure}
	\quad
	\begin{subfigure}[t]{0.475\textwidth}
		\includegraphics[width=\textwidth]{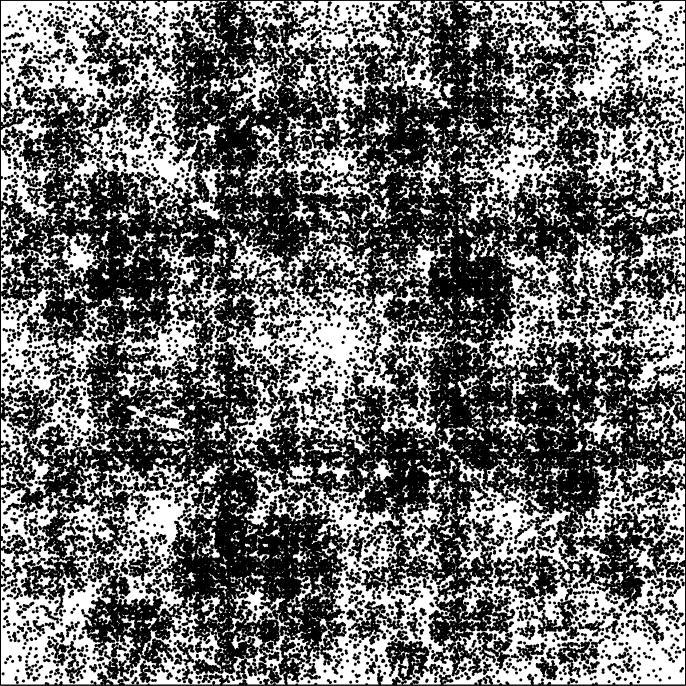}
		\caption{All prime numbers between 7 and $1000000$ taken~$\bmod\ 10$ where the vertices are 1, 3, 7, 9 clockwise from the bottom left}
		\label{fig:cgr_prime_mod10_square}
	\end{subfigure}
	\\
	\begin{subfigure}[t]{0.475\textwidth}
		\includegraphics[width=\textwidth]{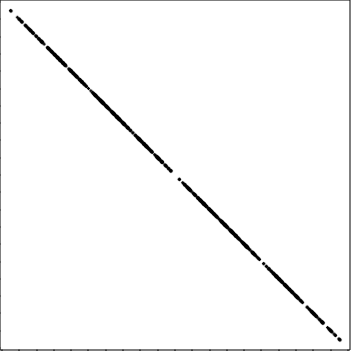}
		\caption{$103 + k$th prime mod 4}
		\label{fig:cgr_ithprime}
	\end{subfigure}
	\quad
	\begin{subfigure}[t]{0.475\textwidth}
		\includegraphics[width=\textwidth]{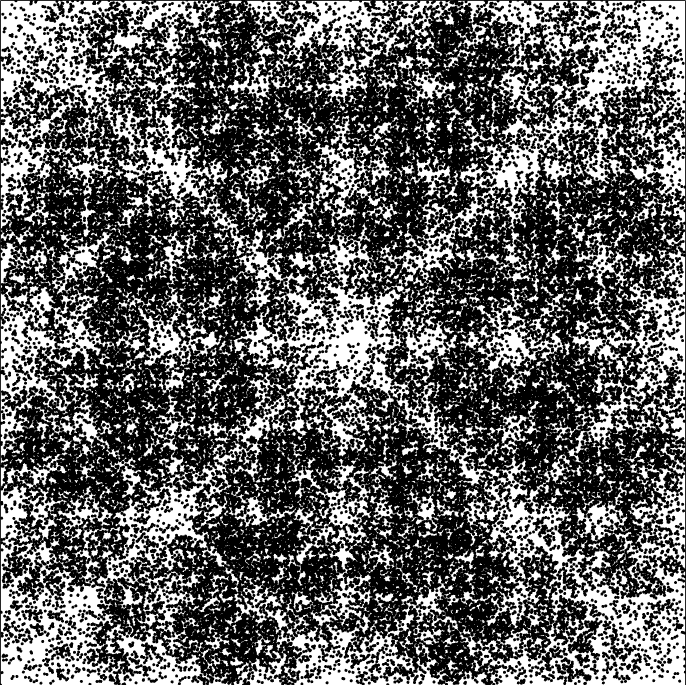}
		\caption{All prime numbers between 7 and $1000000$ taken~$\bmod\ 8$ where the vertices are 1, 3, 5, 7 clockwise from the bottom left}
		\label{fig:cgr_prime_mod_8}
	\end{subfigure}
	\caption{Examples of CGRs of prime numbers}
\end{figure}

We first visualized the prime numbers in a similar way as Figures \ref{fig:digits_pi_mod10} and \ref{fig:cgr_Fibonacci} by taking the sequence of prime numbers modulo 10. This is shown in Figure \ref{fig:cgr_prime_mod_10}. We notice that almost all of the prime numbers ends with a 1, 3, 7, or 9. Because of this observation, we decided to create a square CGR labelling the vertices with these values. The result is shown in Figure \ref{fig:cgr_prime_mod10_square}.

Figure \ref{fig:cgr_ithprime} shows the third visualization of prime numbers: a CGR of the $103+k$th prime numbers $\bmod\ 4$, a diagonal line which spans from the ``1'' corner to the ``3'' corner. As we all know, prime numbers larger than 2 will never be even, so the corner points ``0'' and ``2'' which are representative ``even'' values would never occur, thus creating a straight diagonal line. If we had reassigned the values the coordinate points represent, the plot would have turned out differently; instead of a diagonal line, it could possibly be a horizontal (spanning from coordinate points $(-1, -1)$ to $(1, -1)$ or from $(-1, 1)$ to $(1, 1)$) or vertical line (spanning from $(-1, 1)$ to $(-1, -1)$ or $(1, 1)$ to $(1, -1)$), or a diagonal line the spanning the other two vertices. One needs to be mindful of how the coordinate points are assigned.

Lastly, we can visualize the prime numbers in a different way, as seen in Figure \ref{fig:cgr_prime_mod_8}. Here, we took the sequence of prime numbers less than one million and took $\bmod\ 8$ of the numbers. We thought that this would be interesting because taking any prime number larger than $2\bmod{8}$ will only result in the following numbers: 1, 3, 5, and 7. Because of this, we relabelled the vertices to ``1", ``3", ``5", and ``7" starting from the bottom left corner, going clockwise (i.e.~the coordinate point $(-1, -1)$ is labelled ``1", $(-1, 1)$ is labelled ``3", etc.), and plotted the points according to the chaos game ($r = 0.5$). Surprisingly, this gives us a pattern, apparently meaning that there is a pattern in the sequence of prime numbers~$\bmod\ 8$. However, compare this to Figure \ref{fig:cgr_square}, which, though random, showed patterns. This ``evidence of pattern" is not conclusive! This is possibly simply a Procrustean effect, and we'd have to know more about the distribution of primes to say more.

\subsection{Partial quotients of continued fractions}
We then explored visualizing the sequences taken from the partial quotients of continued fractions, which were taught at the beginning of the course. Continued fractions are fractions written in the following general form~\cite{olds1963continued, corless1992continued}:
\begin{equation}
	a_0 + \cfrac{1}{a_{1} + \cfrac{1}{a_{2} + \cfrac{1}{a_{3} + \cfrac{1}{a_{4} + \cfrac{1}{\ddots}}}}} \>,
\end{equation}
which produces a sequence of numbers, $a_{0} + \left[ a_{1}, a_{2}, a_{3}, a_{4}, \cdots\right]$, called the partial quotients of the continued fraction. This can be demonstrated with an example: we can rewrite \sfrac{9}{7} in the form
\begin{equation}
	\dfrac{9}{7} = 1 + \dfrac{2}{7} = 1 + \dfrac{1}{\sfrac{7}{2}} = 1 + \cfrac{1}{3 + \sfrac{1}{2}} = 1 + \dfrac{1}{3 + \cfrac{1}{1 + \sfrac{1}{1}}} \>.
	\label{eqn:cfrac}
\end{equation}
Here, the partial quotients of $\sfrac{9}{7}$ are the elements of $1 + \left[3, 1, 1\right]$. The students then looked at the continued fractions of other numbers, such as $\sqrt{2}$, $e$ and $\pi$, where the sequences of their partial quotients were used in chaos game representations.

\begin{figure}[ht]
	\centering
	\begin{subfigure}[t]{0.3\textwidth}
		\includegraphics[width=\textwidth]{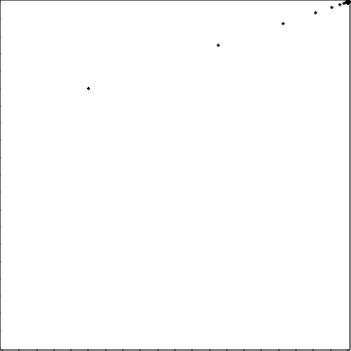}
		\caption{Partial quotients of continued fraction of $\sqrt{2}$ mod 4}
		\label{fig:cgr_cf_sqrt_2}
	\end{subfigure}
	\:
	\begin{subfigure}[t]{0.3\textwidth}
		\includegraphics[width=\textwidth]{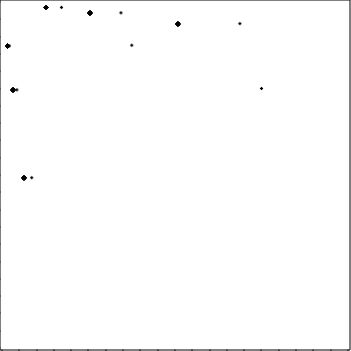}
		\caption{Partial quotients of continued fraction of $e$ mod 4}
		\label{fig:cgr_exp}
	\end{subfigure}
	\:
	\begin{subfigure}[t]{0.3\textwidth}
		\includegraphics[width=\textwidth]{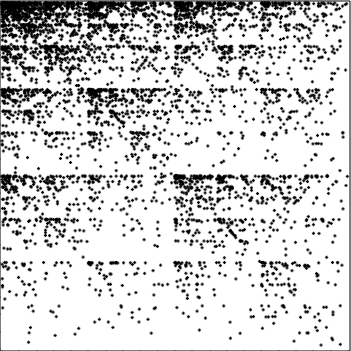}
		\caption{Partial quotients of continued fraction of $\pi$ mod 4}
		\label{fig:cgr_cf_pi}
	\end{subfigure}
	\caption{Examples of CGR for partial quotients of continued fractions of $\sqrt{2}$ (\ref{fig:cgr_cf_sqrt_2}), $e$ (\ref{fig:cgr_exp}), and $\pi$ (\ref{fig:cgr_cf_pi}).}
\end{figure}

One recurring theme of the experimental mathematics course was $\sqrt{2}$, so with this in mind, some students jumped on the opportunity to plot the chaos game representation of the quotients of the continued fraction of $\sqrt{2}$. As we had seen earlier in the course, the sequence goes like
\begin{equation*}
	1 + [2, 2, 2, 2, \ldots, 2] \>.
\end{equation*}
Because all elements of the sequence (apart from the first one) are 2's, it is not surprising to see a (faint) diagonal line with most of the points in the upper right corner, shown in figure \ref{fig:cgr_cf_sqrt_2}. The students also tried the partial quotients of $e$, which is equal to
\begin{equation}
	2 + \left[1, 2, 1, 1, 4, 1, 1, 6, \cdots \right]
\end{equation}
\cite{olds1970simple}, shown in figure \ref{fig:cgr_exp}. What is seen here makes sense as the sequence mostly contains 1's and these alternates with even values, so in this case, either 0 or 2. Therefore, it is clear that there are no points in the lower right corner since that represents the value 3. Unfortunately, these two plots do not look all that impressive; in fact, it is pretty underwhelming. Disappointed with this result, the students decided to experiment with other sequences in which all four coordinates occur.

The students then thought, ``Why not take the partial quotients of the continued fraction of $\pi$? The results \emph{must} be random.'' As shown in figure \ref{fig:cgr_cf_pi}, unexpectedly, there is indeed a pattern, which shows that the sequence of partial quotients of the continued fraction of $\pi$ is not as random as we thought at the beginning of the course. This showed the students that the distribution is not uniform.

In fact, the distribution of partial quotients is very well known, owing to the results of Khinchin~\cite{khinchin1964continued}. As a very startling aside, recently, Bill Gosper has remarked on an amazing identity
\begin{equation}
	\prod k!^{\Delta_{3}\ln k} = K^{\ln 2} \>.,
\end{equation}
where $K$ is Khinchin's constant. We did not inform our students of Khinchin's remarkable results because CGR was only done at the end of the course. Next time, perhaps!

\section{Random vs Pseudorandom\label{sec:random}}
The theory of probability arose much later in the development of science than the theory of dynamical systems and of exact trajectories, and therefore it may be supposed to be more difficult to grasp.  There are many different works on the fundamentals of probability which discuss this difficulty; see for instance~\cite{hamming1991art}, but also~\cite{kac1959probability}. We believe that combining deterministic dynamics with randomness is even more difficult; and that this, in fact, is the main business nowadays of the applied mathematician.  There have been, of course, millions of words written on the topic.  We believe that it is crucial that entering students see some of the discussions of the fundamentals; they need to have a chance to grasp the deeper, most practical aspects of the theory.  We believe that this module offers the instructor a chance to begin those discussions.

One important early part of the discussion is whether or not what the computer produces is ``really'' random.  We alluded earlier to the fact that the poor quality of some early bad random number generators was detected by patterns arising in two-dimensional pictures.  How are the patterns arising in the pictures in this paper different?  Would they arise if ``really'' random numbers were used instead of the pseudorandom numbers generated by computer? [The answer is yes.]  And what is the difference, anyway?  There are knots here that the students (and professors) can tie themselves in: once a sequence of numbers is written down, however randomly it was generated, how can it be random, any more?  It's now perfectly predictable!  The most comforting words that we know about this come from Kolmogorov himself (quoted in full at the end of this paper) which we paraphrase as ``it's only a model.'' Indeed see~\cite{kolmogorov1983logical} for a brief discussion of a \emph{mathematical} foundation for probability.  This applies whether your sequence is ``really'' random, or only ``pseudo'' random.

These discussions are important, in a situation such as that of this paper, where we are trying to tease out deterministic aspects of \emph{apparently} random sequences (or of sequences such as that of DNA which is surely influenced both by random events (mutation, horizontal gene transfer) and very non-random events (selection).  We have used some ``clearly'' non-random (in some sense) sequences such as the digits of $\pi$ to show that we can detect a signature of randomness there; we have used the same techniques to detect a signature of regularity.  To be convinced, the students must be allowed to discuss these issues at some length. In particular, this may be the first time the students have encountered pseudorandom numbers.

There are several high-quality methods for uniform pseudorandom number generation: the most popular includes the \emph{multiple-recursive method} and families formed by \emph{shift-register methods} \cite{dennis2015princeton}. The multiple-recursive method is an extension of the linear congruential method by replacing the first-order linear recursion by one of higher order. One family in the shift-register method generates uniform pseudorandom numbers by means of linear recurring sequences modulo~$2$. Another uses vector recursions modulo~$2$ of higher order; this family includes the very popular Mersenne twister MT19937~\cite{matsumoto1998mersenne}, which is the default for MATLAB's pseudorandom number generating functions, which includes \texttt{rand} (uniformly distributed random number), \texttt{randn} (normally distributed random numbers), and \texttt{randi} (uniformly distributed pseudorandom integers)~\cite{moler2015random}. Other pseudorandom number generators in MATLAB include SIMD-oriented Fast Mersenne Twister, Combined Multiple Recursive, Multiplicative Lagged Fibonacci and Legacy MATLAB generators. Users can use MATLAB's \texttt{rng} function to control which generator to use, along with the seed for the pseudorandom number generator to produce a predictable sequence of numbers.

As mentioned above, one of the pseudorandom number generators MATLAB's \texttt{rand} function uses is the Mersenne twister MT19937 to generate uniform pseudorandom numbers. The MT19937 produces sequences of uniform pseudorandom numbers with period length $2^{19937} - 1$ that possess 623-dimensional equidistribution up to 32 bits accuracy~\cite{dennis2015princeton}. This would seem to create a better ``random" outcome in comparison to an actual sequence of die rolls according to some measures. However, \url{random.org} claims that their ``dice roller'' is, for many purposes, better than pseudorandom number algorithms typically used in computer programs. The randomness in their program comes from atmospheric noise.  

Using a pseudorandom number generator differs to using dice or the atmospheric-sourced generators at \url{random.org}. One gets different results due to the different methods that are used. 
Now, physically-based systems that we believe are random are sometimes more problematic than we think.
For instance, rolling a die is not quite as random as one would expect. Stein's article \cite{stein2012dice} in Inside Science states that dice rolls are not completely random: the initial position of the (fair) die affects the outcome of the die roll. 
And is anything in the atmosphere truly random?  Or is this just a statement of our ignorance?  Kolmogorov's words---it's only a model---are comforting here.  But they help with the use of pseudorandom numbers, as well.

In practice, pseudorandom number generators are indispensable for working with probability, and an enormous amount of work has gone into making them of very high-quality indeed.  We know of no way to distinguish---in practice---results from these good generators from results given by physically-based generators that nearly everyone believes are ``really'' random.

\section{Concluding Remarks}
This module can be used to teach students about a useful visualization of sequences, in which case the emphasis from the instructor could be on interpreting the results, perhaps by using structural similarity or distance maps. The module can also be used to teach elementary programming techniques, in which case the instructor can emphasize programming tools, correctness, efficiency, style, or analysis.

Using elementary sequences of integers makes the module accessible even to first-year students, avoiding difficult biological details: but even so, students are exposed to the deep concepts of pattern and randomness more or less straight away. We have only begun a discussion of what it means to be random, and what it means to have a pattern. Randomness in a model can be a way to hide our ignorance, or it can be a profound statement of our understanding. Modelling viral evolution using randomness~\cite{wahl2017prophage}, where mutations occur in large populations, is clearly warranted; similarly for the vast lengths of DNA sequences. For prime numbers, in one sense clearly not: primes cannot be random, even though they behave in some ways as if they are.
\begin{displayquote}
	``Such considerations may be repeated as often as we like, but it is clear that this procedure will never allow us to be free of the necessity, at the last stage, of referring to probabilities in the primitive imprecise sense of this term.

	It would be quite wrong to think that difficulties of this kind are peculiar in some way to the theory of probability. In the mathematical investigation of actual events, we always make a model of them. The discrepancies between the actual course of events and the theoretical model can, in its turn, be made the subject of mathematical investigation. But for these discrepancies we must construct a model that we will use without formal mathematical analysis of the discrepancies which again would arise in it in actual experiment." \vbox{\hfil--- Andrey Kolmogorov~\cite{aleksandrov1999mathematics}}
\end{displayquote}

\section*{Acknowledgments}
We would like to thank Lila Kari for introducing us to the chaos game representation. We also like to thank the students of the second run of the senior experimental mathematics course for their input and ideas on this topic. We would like to acknowledge financial support from The University of Western Ontario (aka Western University), NSERC, and OGS.

\bibliographystyle{siamplain}
\bibliography{references}

\begin{thebibliography}{10}

\bibitem{aleksandrov1999mathematics}
{\sc A.~D. Aleksandrov, M.~A. Lavrent'ev, et~al.}, {\em Mathematics: Its
  content, methods and meaning}, Courier Corporation, 1999.

\bibitem{almeida2009biological}
{\sc J.~S. Almeida and S.~Vinga}, {\em Biological sequences as pictures--a
  generic two dimensional solution for iterated maps}, BMC bioinformatics, 10
  (2009), p.~100.

\bibitem{barnsley2014fractals}
{\sc M.~F. Barnsley}, {\em Fractals everywhere}, Academic press, 2014.

\bibitem{basu1997chaos}
{\sc S.~Basu, A.~Pan, C.~Dutta, and J.~Das}, {\em Chaos game representation of
  proteins}, Journal of Molecular Graphics and Modelling, 15 (1997),
  pp.~279--289.

\bibitem{borg2005modern}
{\sc I.~Borg and P.~J.~F. Groenen}, {\em Modern multidimensional scaling:
  Theory and applications}, Springer Science \& Business Media, 2005.

\bibitem{chan2018random}
{\sc E.~Y.~S. Chan and R.~M. Corless}, {\em A random walk through experimental
  mathematics}, in Springer Proceedings in Mathematics {\&} Statistics,
  Springer International Publishing, 2020, pp.~203--226,
  \url{https://doi.org/10.1007/978-3-030-36568-4_14},
  \url{https://doi.org/10.1007/978-3-030-36568-4_14}.

\bibitem{corless1992continued}
{\sc R.~M. Corless}, {\em Continued fractions and chaos}, The American
  mathematical monthly, 99 (1992), pp.~203--215.

\bibitem{dennis2015princeton}
{\sc M.~R. Dennis, P.~Glendinning, P.~A. Martin, F.~Santosa, and J.~Tanner},
  {\em The Princeton companion to applied mathematics}, Princeton University
  Press, 2015.

\bibitem{feldman2012chaos}
{\sc D.~P. Feldman}, {\em Chaos and fractals: an elementary introduction},
  Oxford University Press, 2012.

\bibitem{fiser1994chaos}
{\sc A.~Fiser, G.~E. Tusnady, and I.~Simon}, {\em Chaos game representation of
  protein structures}, Journal of molecular graphics, 12 (1994), pp.~302--304.

\bibitem{fisher2012science}
{\sc Y.~Fisher, M.~McGuire, R.~F. Voss, M.~F. Barnsley, R.~L. Devaney, and
  B.~B. Mandelbrot}, {\em The science of fractal images}, Springer Science \&
  Business Media, 2012.

\bibitem{goldman1993nucleotide}
{\sc N.~Goldman}, {\em Nucleotide, dinucleotide and trinucleotide frequencies
  explain patterns observed in chaos game representations of {DNA} sequences},
  Nucleic Acids Research, 21 (1993), pp.~2487--2491.

\bibitem{hamming1991art}
{\sc R.~W. Hamming}, {\em Art of Probability}, Addison Wesley Publishing
  Company, 1991.

\bibitem{oeis}
{\sc O.~F. Inc.}, {\em The on-line encyclopedia of integer sequences}.
\newblock \url{https://oeis.org/}, 2019.

\bibitem{jeffrey1990chaos}
{\sc H.~J. Jeffrey}, {\em Chaos game representation of gene structure}, Nucleic
  Acids Research, 18 (1990), pp.~2163--2170.

\bibitem{jeffrey1992chaos}
{\sc H.~J. Jeffrey}, {\em Chaos game visualization of sequences}, Computers \&
  Graphics, 16 (1992), pp.~25--33.

\bibitem{kac1959probability}
{\sc M.~Kac}, {\em Probability and related topics in physical sciences},
  vol.~1, AMS, 1959.

\bibitem{karamichalis2015investigation}
{\sc R.~Karamichalis, L.~Kari, S.~Konstantinidis, and S.~Kopecki}, {\em An
  investigation into inter-and intragenomic variations of graphic genomic
  signatures}, BMC bioinformatics, 16 (2015), p.~246.

\bibitem{kari2013map}
{\sc L.~Kari, K.~A. Hill, A.~S. Sayem, N.~Bryans, K.~Davis, and N.~S. Dattani},
  {\em Map of life: Measuring and visualizing species' relatedness with"
  molecular distance maps"}, arXiv preprint arXiv:1307.3755,  (2013).

\bibitem{kari2015mapping}
{\sc L.~Kari, K.~A. Hill, A.~S. Sayem, R.~Karamichalis, N.~Bryans, K.~Davis,
  and N.~S. Dattani}, {\em Mapping the space of genomic signatures}, PloS one,
  10 (2015), p.~e0119815.

\bibitem{khinchin1964continued}
{\sc A.~Khinchin and T.~Teichmann}, {\em Continued fractions}, Physics Today,
  17 (1964), p.~70.

\bibitem{knuth1997seminumerical}
{\sc D.~E. Knuth}, {\em Seminumerical Algorithms}, vol.~2 of The Art of
  Computer Programming, Addison Wesley, Reading, MA, 3~ed., 1997.

\bibitem{kolmogorov1983logical}
{\sc A.~N. Kolmogorov}, {\em On logical foundations of probability theory}, in
  Probability theory and mathematical statistics, Springer, 1983, pp.~1--5.

\bibitem{matsumoto1998mersenne}
{\sc M.~Matsumoto and T.~Nishimura}, {\em Mersenne twister: a 623-dimensionally
  equidistributed uniform pseudo-random number generator}, ACM Transactions on
  Modeling and Computer Simulation (TOMACS), 8 (1998), pp.~3--30.

\bibitem{moler2015random}
{\sc C.~B. Moler}, {\em Random number generators, {M}ersenne {T}wiser}.
\newblock In \textit{{C}leve's {C}orner: {C}leve {M}oler on Mathematics and
  Computing}. Available at
  \url{https://blogs.mathworks.com/cleve/2015/04/17/random-number-generator-mersenne-twister/},
  2015.

\bibitem{olds1963continued}
{\sc C.~D. Olds}, {\em Continued fractions}, vol.~18, Random House New York,
  1963.

\bibitem{olds1970simple}
{\sc C.~D. Olds}, {\em The simple continued fraction expansion of e}, American
  Mathematical Monthly,  (1970), pp.~968--974.

\bibitem{simon1991calculation}
{\sc I.~Simon, L.~Glasser, and H.~A. Scheraga}, {\em Calculation of protein
  conformation as an assembly of stable overlapping segments: application to
  bovine pancreatic trypsin inhibitor.}, Proceedings of the National Academy of
  Sciences, 88 (1991), pp.~3661--3665.

\bibitem{stein2012dice}
{\sc B.~P. Stein}, {\em Dice rolls are not completely random}.
\newblock Available at
  \url{https://www.insidescience.org/news/dice-rolls-are-not-completely-random}
  (September 22, 2012).

\bibitem{wahl2017prophage}
{\sc L.~M. Wahl and T.~Pattenden}, {\em Prophage provide a safe haven for
  adaptive exploration in temperate viruses}, Genetics, 206 (2017),
  pp.~407--416.

\bibitem{wang2004image}
{\sc Z.~Wang, A.~C. Bovik, H.~R. Sheikh, and E.~P. Simoncelli}, {\em Image
  quality assessment: from error visibility to structural similarity}, IEEE
  trans. image process., 13 (2004), pp.~600--612.

\bibitem{wickelmaier2003introduction}
{\sc F.~Wickelmaier}, {\em An introduction to {MDS}}, Sound Quality Research
  Unit, Aalborg University, Denmark, 46 (2003).

\end{thebibliography}
\end{document}